\documentclass[12pt]{article}
\usepackage{amsfonts,amsmath,amssymb,amsthm}
\usepackage{color}
\usepackage{graphicx}
\usepackage[tight]{subfigure}
\usepackage{epic,eepic}
\usepackage{fullpage,multicol}

\usepackage{times}

\newcommand{\aff}{a_{\rm ff}}
\newcommand{\afb}{a_{\rm fb}}
\newcommand{\eps}{\varepsilon}
\newcommand{\lmax}{\Lambda_{\max}}
\newcommand{\bA}{{\bf A}}
\newcommand{\bT}{{\bf T}}

\newcommand{\heading}[1]{\subsection*{\normalsize #1}}

\newtheorem{theorem}{Theorem}

\newcommand{\kaption}[1]{\caption{\small #1}}

\title{Shear-Induced Chaos}

\author{Kevin K. Lin\thanks{E-mail: {\bf klin@cims.nyu.edu}.
    K.~L. is supported by an NSF postdoctoral fellowship.}
\ and Lai-Sang Young\thanks{E-mail: {\bf lsy@cims.nyu.edu}.
    L.-S.~Y. is supported by a grant from the NSF.}\\ Courant
  Institute of Mathematical Sciences\\ New York University}

\date{May 22, 2007}

\begin{document}
\maketitle
\begin{abstract}
Guided by a geometric understanding developed in earlier works
of Wang and Young, we carry out some numerical studies of
shear-induced chaos. The settings considered include periodic
kicking of limit cycles, random kicks at Poisson times, and
continuous-time driving by white noise. The forcing of a
quasi-periodic model describing two coupled oscillators is also
investigated.  In all cases, positive Lyapunov exponents are
found in suitable parameter ranges when the forcing is suitably
directed.
\end{abstract}

%%--------------------------------%%

\section*{Introduction}

This paper presents a series of numerical studies which
investigate the use of shear in the production of chaos.  The
phenomenon in question can be described roughly as follows: An
external force is applied to a system with tame, nonchaotic
dynamics.  If the forcing is strategically applied to interact
with the shearing in the underlying dynamics, it can sometimes
lead to the folding of phase space, which can in turn lead to
positive Lyapunov exponents for a large set of initial
conditions.  This phenomenon, which we call {\em shear-induced
  chaos}, occurs in a wide variety of settings, including
periodically-forced oscillators.  For a topic as general as
this, it is difficult to compile a reasonable set of references.
We have not attempted to do that, but mention that the first
known observation of a form of this phenomenon was by van der
Pol and van der Mark 80 years ago {\cite{van der pol}}.  Other
references related to our work will be mentioned as we go along.

The starting point of the present work is a series of papers by
Wang and Young {\cite{WY1,WY2,WY3,WY4}}.  In these papers, the
authors devised a method for proving the existence of strange
attractors and applied their techniques to some natural
settings; two of their main examples are periodically-kicked
oscillators and systems undergoing Hopf bifurcations.  They
identified a simple geometric mechanism to explain how the
chaotic behavior comes about.  Because of the perturbative
nature of their analysis, however, the kicks in these results
have to be followed by very long periods of relaxation.  In
other words, the chaos in these attractors develop on a very
slow time scale.  Relevant parts of the results of Wang and
Young are reviewed in Sect.~\ref{wang-young}.

The aim of the present paper is to study shear-induced chaos in
situations not accessible by current analytic tools.  We believe
that this phenomenon is widespread, meaning it occurs for large
sets of parameters, and that it is robust, meaning it does not
depend sensitively on the type of forcing or even background
dynamics as long as certain geometric conditions are met.  We
validate these ideas through a series of numerical studies in
which suitable parameters are systematically identified
following ideas from {\cite{WY2}} and {\cite{WY3}}.  Four
separate studies are described in
Sects.~\ref{study1}--~\ref{study4}.  The first three studies
involve an oscillator driven by different types of forcing (both
deterministic and stochastic); in these studies, the unforced
system is a simple linear shear flow model.  In the fourth
study, the unforced dynamics are that of a coupled oscillator
system described by a (periodic or quasi-periodic) flow on the
$2$-torus.

The linear shear flow used in Studies 1--3 has been studied
independently in {\cite{zaslav}} and {\cite{WY2}}.  It is the
simplest system known to us that captures all the essential
features of typical oscillator models relevant to shear-induced
chaos.  Moreover, these features appear in the system in a way
that is easy to control, and the effects of varying each are
easy to separate.  This facilitates the interpretation of our
theoretical findings in more general settings in spite of the
fact that numerical studies necessarily involve specific models.

We mention that our results on shear flows are potentially
applicable to a setting not discussed here, namely that of the
advection and mixing of passive scalar tracers in (weakly
compressible) flows.

Finally, we remark that this work exploits the interplay between
deterministic and stochastic dynamics in the following way: The
geometry in deterministic models are generally more clear-cut.
It enables us to extract more readily the relationship between
quantities and to deduce the type of results these relationships
may lead.  Results for stochastic models, on the other hand,
tend to be more {\em provable} than their counterparts in
deterministic models, where competing scenarios lead to very
delicate dependences on parameters.  Our numerical results on
stochastic forcing in Studies 2--4 point clearly to the
possibility of (rigorous) theorems, some versions of which, we
hope, will be proved in the not too distant future.

%%--------------------------------%%

\section{Rigorous Results and Geometric Mechanism}
\label{wang-young}

In this section, we review some rigorous results of Wang and
Young (mainly \cite{WY2,WY3}, also \cite{WY1,WY4}) and the
geometric mechanism for producing chaos identified in the first
two of these papers. We will focus on the case of limit cycles,
leaving the slightly more delicate case of supercritical Hopf
bifurcations to the reader. The material summarized in this
section form the starting point for the numerical investigations
in the present paper.

\subsection{Strange Attractors from Periodically-Kicked Limit Cycles}
\label{strangeattractors}

Consider a smooth flow $\Phi_t$ on a finite dimensional
Riemannian manifold $M$ (which can be ${\mathbb R}^d$), and let
$\gamma$ be a {\em hyperbolic limit cycle}, {\it i.e.} $\gamma$
is a periodic orbit of $\Phi_t$ with the property that if we
linearize the flow along $\gamma$, all of the eigenvalues
associated with directions transverse to $\gamma$ have strictly
negative real parts.  The {\em basin of attraction} of $\gamma$,
${\cal B}(\gamma)$, is the set $\{x \in M
:\Phi_t(x)\to\gamma\mbox{ as $t\to\infty$}\}$.  It is well known
that hyperbolic limit cycles are robust, meaning small
perturbations of the flow will not change its dynamical picture
qualitatively.

A {\em periodically-kicked oscillator} is a system in which
``kicks" are applied at periodic time intervals to a flow $\Phi_t$
with a hyperbolic limit cycle. For now let us think of a ``kick"
as a mapping $\kappa: M \to M$. If kicks are applied $T$ units
of time apart, then the time evolution of the kicked system 
can be captured by iterating its time-$T$ map
$F_T = \Phi_T \circ \kappa$. If there is a
neighborhood $\cal U$ of $\gamma$ such that $\kappa({\cal U})
\subset {\cal B}(\gamma)$, and the relaxation time is long
enough that points in $\kappa({\cal U})$ return to $\cal U$,
{\it i.e.}, $F_T({\cal U}) \subset {\cal U}$,
then $\Gamma = \cap_{n \geq 0} F_T^n({\cal U})$ is an attractor
for the periodically kicked system $F_T$. In a sense,
$\Gamma=\Gamma(\kappa, T)$ is what becomes of the limit cycle
$\gamma$ when the oscillator is periodically kicked.  Since
hyperbolic limit cycles are robust, $\Gamma$ is a slightly
perturbed copy of $\gamma$ if the kicks are weak.  We call it an
``invariant circle.''  Stronger kicks may ``break" the invariant circle,
leading to a more complicated invariant set. 
Of interest in this paper is when $\Gamma$ is a
strange attractor, {\it i.e.}, when the dynamics in $\cal U$
exhibit sustained, observable chaos.

Two theorems are stated below.  Theorem 1 is an abstract result,
the purpose of which is to emphasize the generality of the
phenomenon.  Theorem 2 discusses a concrete situation intended
to make transparent the relevance of certain quantities.  Let Leb$(\cdot)$
denote the Lebesgue measure of a set.

\begin{theorem}  {\rm\cite{WY3}}
Let $\Phi_t$ be a $C^4$ flow with a hyperbolic limit cycle
$\gamma$. Then there is an open set of kick maps $\cal K$ with
the following properties: For each $\kappa \in {\cal K}$, there
is a set $\Delta = \Delta(\kappa) \subset {\mathbb R}^+$ with
{\rm Leb}$(\Delta)>0$ such that for each $T \in \Delta$,
$\Gamma$ is a ``strange attractor'' of $F_T$.
\end{theorem}

The term ``strange attractor'' in the theorem has a well-defined
mathematical meaning, which we will discuss shortly.  But first
let us take note of the fact that this result applies to all systems
with hyperbolic limit cycles, independent of dimension or other
specifics of the defining equations.  Second, we remark that
the kicks in this theorem are very infrequent, {\it i.e.} $T \gg
1$, and that beyond a certain $T_0$, the set $\Delta$ is roughly
periodic with the same period as the cycle $\gamma$.

The term ``strange attractor'' in Theorem 1 is used as
short-hand for an attractor with a package of well defined
dynamical properties.  These properties were established for a
class of rank-one attractors (see {\cite{WY1}} for the
2-dimensional case; a preprint for the $n$-dimensional case will
appear shortly).  In {\cite{WY1, WY4}}, the authors identified a
set of conditions that implies the existence of such attractors,
and the verification of the conditions in {\cite{WY1, WY4}} in
the context of Theorem 1 is carried out in {\cite{WY3}} (see
also {\cite{gwy}} and {\cite{oksa+wang}} for other applications
of these ideas).  We refer the reader to the cited papers for
more details, and mention only the following three
characteristics implied by the term ``strange attractor" in this
section.
\begin{itemize}

\item[(1)] There is a set $\cal V$ of full Lebesgue measure in
  the basin of attraction of $\Gamma$ such that orbits starting
  from every $x \in \cal V$ have (strictly) positive Lyapunov
  exponents.

\item[(2)] $F_T$ has an ergodic SRB measure $\mu$, and 
for every continuous observable $\varphi$, 
$$\frac{1}{n}\sum_{i=0}^{n-1} \varphi(F_T^i (x)) \to \int
\varphi\ d \mu \qquad\mbox{as $n \to \infty$ for every $x \in
  \cal V$.}
$$

\item[(3)] The system $(F_T, \mu)$ is mixing; in fact, it has
  exponential decay of correlations for H\"older continuous
  observables.

\end{itemize}

An important remark before leaving Theorem 1: Notice that the
existence of ``strange attractors'' is asserted for $F_T$ for
only a positive measure set of $T$, not for all sufficiently
large $T$.  This is more a reflection of reality than a weakness
of the result. For large enough $T$ in the complement of
$\Delta$, the attractor is guaranteed to contain horseshoes, the
presence of which will lead to some semblance of chaotic
behavior. After a transient, however, typical orbits may (or may
not) tend to a stable equilibrium. If they do, we say $F_T$ has
{\it transient chaos}. This is to be contrasted with properties
(1)--(3) above, which represent a much stronger form of chaos.

The next result has an obvious analog in $n$-dimensions (see
\cite{WY3}), but the 2-D version illustrates the point.

\begin{theorem}  {\rm\cite{WY3}}
Consider the system
\begin{equation}
\begin{array}{rcl}
\dot{\theta} & =& 1 + \sigma y\\
\dot{y} &=& -\lambda y + A\cdot H(\theta) \cdot\sum_{n=0}^\infty \delta(t-nT)
\end{array}
\label{thm2}
\end{equation}
where $(\theta, y) \in S^1 \times {\mathbb R}$ are coordinates
in the phase space, $\lambda, \sigma, A >0$ are constants, and
$H: S^1 \to {\mathbb R}$ is a {\rm nonconstant} smooth function.
If the quantity
$$
\frac{\sigma}{\lambda} \cdot A \ \equiv \
\frac{\mbox{\rm shear}}{\mbox{\rm contraction rate}} \cdot
\mbox{\rm kick ``amplitude''}
$$
is sufficiently large (how large depends on the forcing function
$H$), then there is a positive measure set $\Delta \subset
{\mathbb R}^+$ such that for all $T \in \Delta$, $F_T$ has a
strange attractor in the sense above.
\end{theorem}

\begin{figure}
\begin{center}
\includegraphics[bb=0 0 307 117]{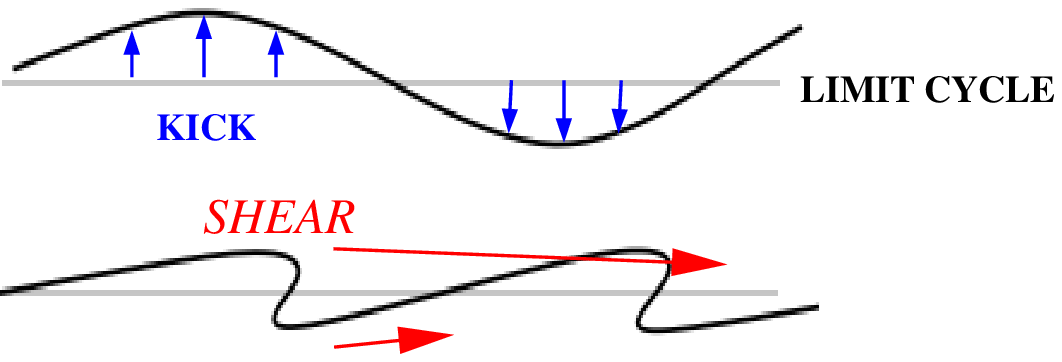}
\end{center}
\kaption{The stretch-and-fold action of a kick followed by
  relaxation in the presence of shear.}
\label{shear}
\end{figure}

Here, the term involving $H(\theta)$ defines the kick, and
$\gamma=S^1\times\{0\}$.  We explain intuitively the
significance of the quantity $\frac{\sigma}{\lambda} A$. As
noted earlier, to create a strange attractor, it is necessary to
``break'' the limit cycle.  The more strongly attractive
$\gamma$ is, the harder it is to break.  From this we see the
advantage of having $\lambda$ small. By the same token, a
stronger forcing, {\it i.e.}, larger $A$, helps.  The role of
$\sigma$, the {\it shear}, is explained pictorially in
Fig.~\ref{shear}: Since the function $H$ is required to be
nonconstant, let us assume the kick drives some points on the
limit cycle $\gamma$ up and some down, as shown.  The fact that
$\sigma$ is positive means that points with larger
$y$-coordinates move faster in the $\theta$-direction. During
the relaxation period, the ``bumps'' created by the kick are
stretched as depicted.  At the same time, the curve is attracted
back to the limit cycle.  Thus, the combination of kicks and
relaxation provides a natural mechanism for repeated stretching
and folding of the limit cycle.  Observe that the larger the
differential in speed in the $\theta$-direction, {\it i.e.} the
larger $\sigma$, and the slower the return to $\gamma$, {\it
  i.e.} the smaller $\lambda$, the more favorable the conditions
are for this stretch-and-fold mechanism.

%%%%%%%%%%%%%%%%%%%%%%%%%%%%%%%%%
\subsection{Geometry and Singular Limits}
\label{sec1.2}

In Eq.~(\ref{thm2}), the quantities $\lambda$, $\sigma$ and $A$
appear naturally. But what about in general limit cycles, where
the direction of the kicks vary? What, for example, will play
the role of $\sigma$, or what we called shear in Eq.~(\ref{thm2})? 
The aim of this subsection is to
shed light on the general geometric picture, and to explain how
the dynamics of $F_T$ for large $T$ can be understood.

\heading{Geometry of $F_T$ and the Strong Stable Foliation}

Let $\gamma$ be a hyperbolic limit cycle as in the beginning of
Sect.~1.1.  Through each $x \in \gamma$ passes the {\it strong
  stable manifold} of $x$, denoted $W^{ss}(x)$ {\cite{gh}}. By
definition, $W^{ss}(x)=\{y \in M: d(\Phi_t(y), \Phi_t(x)) \to 0$
as $t \to \infty \}$; the distance between $\Phi_t(x)$ and
$\Phi_t(y)$ in fact decreases exponentially. Some basic
properties of strong stable manifolds are: (i) $W^{ss}(x)$ is a
codimension one submanifold transversal to $\gamma$ and meets
$\gamma$ at exactly one point, namely $x$; (ii)
$\Phi_t(W^{ss}(x)) = W^{ss}(\Phi_t(x))$, and in particular, if
the period of $\gamma$ is $p$, then
$\Phi_p(W^{ss}(x))=W^{ss}(x)$; and (iii) the collection
$\{W^{ss}(x), x \in \gamma\}$ foliates the basin of attraction
of $\gamma$, that is to say, they partition the basin into
hypersurfaces.

\begin{figure}
\begin{center}
\includegraphics[bb=0 0 379 107]{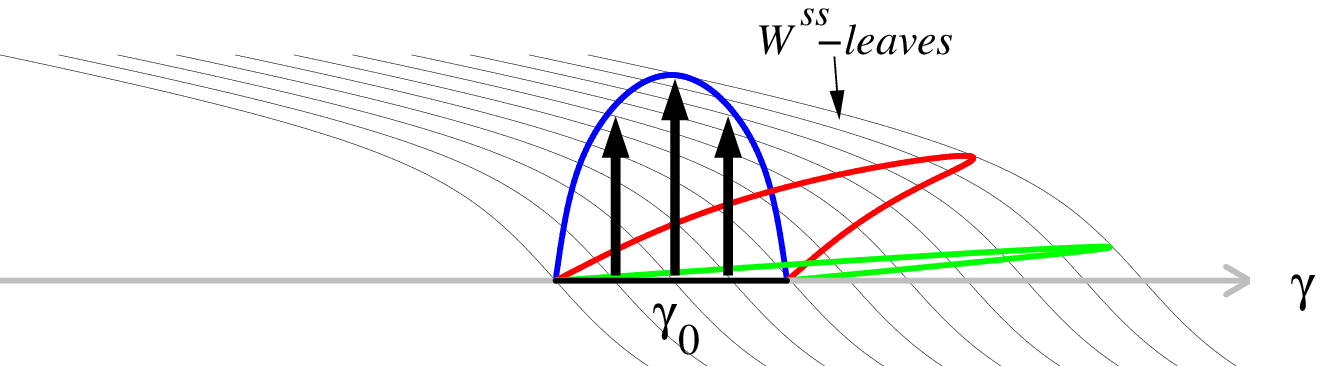}
\end{center}
\kaption{Geometry of folding in relation to the
  $W^{ss}$-foliation.  Shown are the kicked image of a segment
  $\gamma_0$ and two of its subsequent images under
  $\Phi_{np}$.}
\label{wss}
\end{figure}

We examine next the action of the kick map $\kappa$ in relation
to $W^{ss}$-manifolds. Fig.~\ref{wss} is analogous to
Fig.~\ref{shear}; it shows the image of a segment $\gamma_0$ of
$\gamma$ under $F_T = \Phi_T \circ \kappa$. For illustration 
purposes, we assume $\gamma_0$ is kicked upward with
its end points held fixed, and assume $T=np$ for some 
$n \in {\mathbb Z}^+$ (otherwise the
picture is shifted to another part of $\gamma$ but is
qualitatively similar). Since $\Phi_{np}$ leaves each
$W^{ss}$-manifold invariant, we may imagine that during
relaxation, the flow ``slides'' each point of the curve
$\kappa(\gamma_0)$ back toward $\gamma$ along $W^{ss}$-leaves.
In the situation depicted, the effect of the folding is evident.

Fig.~\ref{wss} gives considerable insight into what types of
kicks are conducive to the formation of strange attractors.
Kicks along $W^{ss}$-leaves or in directions roughly parallel to
the $W^{ss}$-leaves will not produce strange attractors, nor
will kicks that essentially carry one $W^{ss}$-leaf to another.
What causes the stretching and folding is the {\it variation} in
how far points $x\in\gamma$ are moved by $\kappa$ as measured in
the direction transverse to the $W^{ss}$-leaves.  Without
attempting to give a more precise characterization, we will
refer to the type of chaos that results from the geometry above
as {\it shear-induced chaos}.  We emphasize that the occurrence
of shear-induced chaos relies on the interplay between the
geometries of the kicks and the underlying dynamical structures.

Returning to the concrete situation of Theorem 2, since
Eq.~(\ref{thm2}) without the kick term is linear, it is easy to
compute strong stable manifolds.  In $(\theta, y)$-coordinates,
they are lines with slope $-\lambda/\sigma$.  Variations in kick
distances here are guaranteed by the fact that $H$ is
nonconstant. With $H$ fixed, it is clear that the larger
$\sigma/\lambda$ and $A$, the greater these variations.  Note
that the use of the word kick ``amplitude" in the statement of
Theorem 2 is a little misleading, for it is not the amplitude of
the kicks {\it per se} that leads to the production of chaos.

%%%%%%%%%%%%
\heading{Singular Limits of $F_T$ as $T \to \infty$}

When $T \gg 1$, {\it i.e.} when kicks are very infrequent, the
map $F_T$ sends a small tube ${\cal U}_T$ around $\gamma$ back
into itself.  This is an example of what is called a {\it
  rank-one map} in \cite{WY4}.  Roughly speaking, a rank-one map
is a smooth map whose derivative at each point is strongly
contractive in all but one of the directions.  Rank-one maps can
be analyzed using perturbative methods if they have well-defined
``singular limits.''  In the context of limit cycles, these
singular limits do exist; they are a one-parameter family of
maps $\{f_a:\gamma\circlearrowleft\}$ obtained by letting
$T\to\infty$ in the following way: For each $a\in[0,p)$ (recall
  that $p =\mbox{ period of $\gamma$}$), let
\begin{equation}
f_a(x) \ := \ \lim_{n \to
  \infty}\Phi_{np+a}(\kappa(x))\qquad\mbox{ for all }
x\in\gamma\ .
\end{equation}
Equivalently, $f_a(x)$ is the unique point $y \in \gamma$ such
that $\kappa(x) \in W^{ss}(y)$.  Notice that $f_a(x) =
f_0(x)+a\mbox{ (mod 1)}$, where we identify $\gamma$ with
$[0,1]$ (with the end points identified).  For Eq.~(\ref{thm2}),
$f_a$ is easily computed to be
\begin{equation}
f_a(\theta)=\theta + a + \frac{\sigma}{\lambda} A \cdot H(\theta),
\label{sing}
\end{equation}
where the right side should again be interpreted as mod 1.  (In
the setting of driven oscillators, singular limits are sometimes
known as ``phase resetting curves''; they have found widespread
use in {\em e.g.}  mathematical biology {\cite{winfree,guck75}}.)

It is shown in \cite{WY1,WY2,WY3,WY4} that a great deal of
information on the attractor $\Gamma$ of $F_T$ for $T \gg 1$ can
be recovered from these singular limit maps.  The results are
summarized below.  These results hold generally, but as we step
through the 3 cases below, it is instructive to keep in mind
Eq.~(\ref{thm2}) and its singular limit (\ref{sing}), with
$\frac{\sigma}{\lambda} A$ increasing as we go along:
\begin{itemize}

\item[(i)] If $f_a$ is injective, {\it i.e.}, it is a circle
  diffeomorphism, the attractor $\Gamma$ for $F_T$ is an
  invariant circle.  This happens when the kicks are aimed in
  directions that are ``unproductive" (see above), or when their
  effects are damped out quickly.  In this case, the competing
  scenarios on $\Gamma$ are quasi-periodicity and ``sinks,''
  {\it i.e.} the largest Lyapunov exponent of $F_T$ is zero or
  negative.

\item[(ii)] When $f_a$ loses its injectivity, the invariant
  circle is ``broken''. When that first happens, the expansion
  of the 1-D map $f_a$ is weak, and all but a finite number of
   trajectories tend to sinks. This translates into a gradient
  type dynamics for $F_T$.

\item[(iii)] If $f_a$ is sufficiently expanding away from its
  critical points, $\Gamma$ contains horseshoes for all large
  $T$. For an open set of these $T$, the chaos is transient,
  while on a positive measure set, $F_T$ has a strange attractor
  with the properties described in Sect. 1.1. These are the two
  known competing scenarios. (They may not account for all $T$.)
  Since $F_T \approx F_{T+np}$ for large $T$, both sets of
  parameters are roughly periodic.

  The analyses in the works cited suggest that when horseshoes
  are first formed, the set of parameters with transient chaos
  is more dominant.  The stronger the expansion of $f_a$, the
  larger the set of parameters with strange attractors. In the
  first case, the largest Lyapunov exponent of $F_T$ may appear
  positive for some time (which can be arbitrarily long) before
  turning negative.  In the second case, it stays positive
  indefinitely.

\end{itemize}

%%%%%%%%%%%%%%%%%%%%%%%%%%%%%%%%%%
\subsection{Limitations of Current Analytic Techniques}

In hyperbolic theory, there is, at the present time, a very large discrepancy 
between what is thought to be true and what can be proved.
Maps that are dominated by stretch-and-fold behavior are generally
thought to have positive Lyapunov exponents -- although this reasoning
is also known to come with the following caveat: Maps whose
derivatives expand in certain directions tend to contract in other directions, 
and unless the expanding and contracting directions are well separated
(such as in Anosov systems), the contractive directions can conspire
to form sinks. This is how the transient chaos described in
Sect. 1.2 comes about.
Still, if the expansion is sufficiently strong, one would
expect that positive Lyapunov exponents are more likely to prevail
-- even though for any one map the outcome can go either way. 
{\it Proving} results of this type is a different matter. 
Few rigorous results exist for systems for which one has 
no {\it a priori} knowledge of invariant cones, and invariant cones are
unlikely in shear-induced chaos.

The rigorous results reviewed in the last two subsections have
the following limitations: (i) They pertain to $F_T$ for only
very large $T$. This is because the authors use a perturbative
theory that leans heavily on the theory of 1-D maps.  No
non-perturbative analytic tools are currently available.  (ii) A
larger than necessary amount of expansion is required of the
singular limit maps $f_a$ in the proof of strange
attractors. This has to do with the difficulty in locating
suitable parameters called Misiurewicz points from which to
perturb. (This problem can be taken care of, however, by
introducing more parameters.) We point out that (i) and (ii)
together exacerbate the problem: $f_a$ is more expanding when
$\lambda$ is small, but if $F_T=\Phi_T \circ \kappa$ is to be
near its singular limit, then $e^{-\lambda T}$ must be very
small, {\it i.e.} $\lambda T$ must be very large.

\medskip
That brings us to the present paper, the purpose of which is to
supply numerical evidence to support some of our conjectured
ideas regarding situations beyond the reach of the rigorous work
reviewed. Our ideas are based on the geometry outlined in
Sect.~\ref{sec1.2}, but are not limited to periodic kicks or to
the folding of limit cycles.

%%%%%%%%%%%%%%%%%%%%%%%%%%%%%%%%%%
%%%%%%%%%%%%%%%%%%%%%%%%%%%%%%%%
\section{Study 1: Periodically-Kicked Oscillators}
\label{study1}

Our first model is the periodic kicking of a linear shear flow 
with a hyperbolic limit cycle. The setting is as in Theorem 2 with 
$H(\theta)=\sin(2\pi\theta)$, {\it i.e.}, we consider
\begin{equation}
\begin{array}{ccl}
\dot{\theta} &=& 1 + \sigma y\ ,\\
\dot{y} &=& -\lambda y + A\cdot \sin(2\pi\theta)
\cdot\sum_{n=0}^\infty \delta(t-nT)\ ,\\
\end{array}
\label{model1}
\end{equation}
where $(\theta, y) \in S^1 \times {\mathbb R}$, $S^1 \equiv
[0,1]$ with the two end points of $[0,1]$ identified. In the
absence of kicks, {\it i.e.}, when $A=0$, $\Phi_t(z)$ tends to
the limit cycle $\gamma=S^1 \times\{0\}$ for all $z \in S^1 \times
{\mathbb R}$.  As before, the attractor in the kicked system is
denoted by $\Gamma$.  The parameters of interest are:
\begin{align*}
\sigma &= \mbox{amount of shear,}\\
\lambda &= \mbox{damping or rate of contraction to $S^1 \times\{0\}$,}\\
A &= \mbox{amplitude of kicks, and}\\
T &= \mbox{time interval between kicks.}
\end{align*}
Our aim here is to demonstrate that the set of parameters with
chaotic behavior is considerably larger than what is guaranteed
by the rigorous results reviewed in Sect.~\ref{wang-young}, and
to gain some insight into this parameter set.  By ``chaotic
behavior,'' we refer in this section to the property that $F_T$
has a positive Lyapunov exponent for orbits starting from a
``large'' set of initial conditions, {\em i.e.}  a set of full
or nearly full Lebesgue measure in the basin of attraction of
$\Gamma$.  More precisely, we {\it assume} that such Lyapunov
exponents are well defined, and proceed to compute the largest
one, which we call $\Lambda_{\max}$.

We begin with some considerations relevant to the search for
parameters with $\Lambda_{\max} > 0$:
\begin{itemize}

\item[(a)] It is prudent, in general, to ensure that orbits do
  not stray too far from $\gamma$. This is because while the basin of
  attraction of $\gamma$ in this model is the entire phase
  space, the basin is bounded in many other situations.  We
  therefore try to keep $\Gamma \subset \{|y|<b\}$ with
  relatively small $b$. This is guaranteed if $A$ is small enough
  that $e^{-\lambda T}(b+A)<b$; the bound is improved
   if, for example, no point
  gets kicked to maximum amplitude two consecutive iterates.

\item[(b)] Let $(\theta_T, y_T)=F_T(\theta_0, y_0)$. A simple
  computation gives
\begin{equation}
\label{time-t map}
\begin{array}{rcl}
\theta_T & = & \theta_0 + T + \frac{\sigma}{\lambda} \cdot 
[ y_0 + A \sin(2\pi\theta_0)] \cdot (1 - e^{-\lambda
  T})\qquad\mbox{(mod 1)}\ ,\\
y_T & = & e^{-\lambda T} [y_0 + A \sin(2\pi\theta_0)]\ .
\end{array}
\end{equation}
For $b$ relatively small, we expect the number $\frac{\sigma
  A}{\lambda}(1 - e^{-\lambda T})$ to be a good indicator of
chaotic behavior: if it is large enough, then $F_T$ folds the
annulus $\{|y|<b\}$ with two turns and maps it into itself.
The larger this number, the larger the folds, meaning 
the more each of the monotonic parts of the image wraps 
around in the $\theta$-direction.
\end{itemize}

\paragraph{Summary of Findings.}
{\em 
\begin{itemize}
\item[(i)] With the choice of parameters guided by (a) and (b) above,
  we find that as soon as the folding described in (b) is 
  definite, $F_T$ becomes ``possibly chaotic'', 
  meaning $\Lambda_{\max}$
  is seen numerically to oscillate (wildly) between positive and
  negative values as $T$ varies. We interpret this to be due to
  competition between transient and sustained chaos; see (iii)
  in Sect. 1.2.  For larger $\frac{\sigma}{\lambda} A$, {\it i.e.},
  as the stretching is stronger, and for $T$ beyond an initial range,
  this oscillation stops and
  $\Lambda_{\max}$ becomes definitively positive for all the
  values of $T$ computed. 
\item[(ii)] As for the range of parameters with chaotic
  dynamics, we find that $\Lambda_{\max}>0$ occurs under fairly
  modest conditions, {\it e.g.}, for $\frac{\sigma}{\lambda} A =
  3$, we find $\Lambda_{\max}>0$ starting from about $T \approx
  3$, which is very far from the ``$T\to \infty$'' in rigorous
  proofs. Also, while shear-induced chaos is often associated
  with weak damping, we find that the phenomenon occurs as well
  for larger $\lambda$, {\it e.g.}, for $\lambda \sim 1$,
  provided its relation to the other parameters are favorable.

\end{itemize}
}

\begin{figure}
\begin{center}
\includegraphics[bb=0in 0in 4in 2.5in]{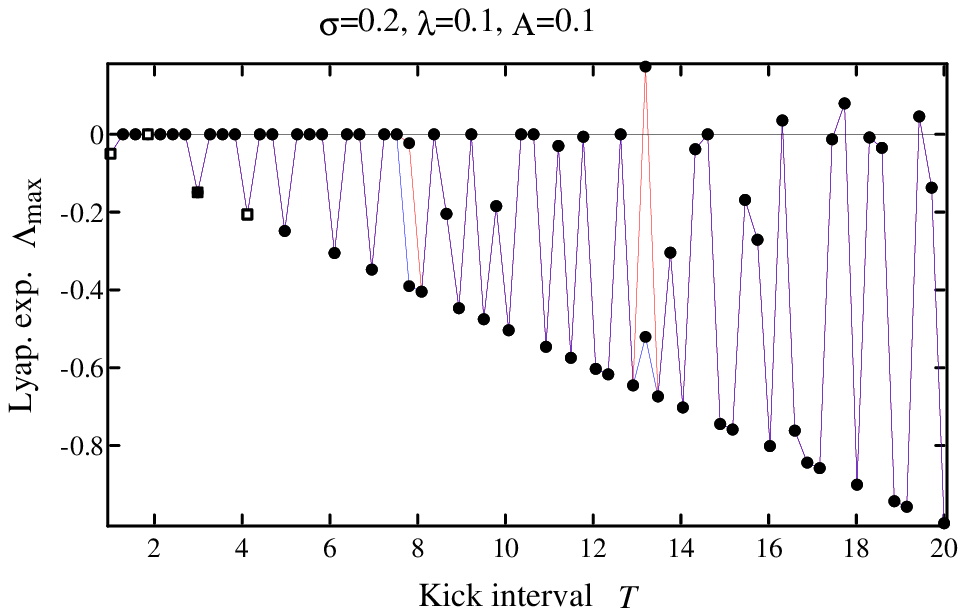}\vspace{18pt}\\
\includegraphics[bb=0in 0in 4in 2.5in]{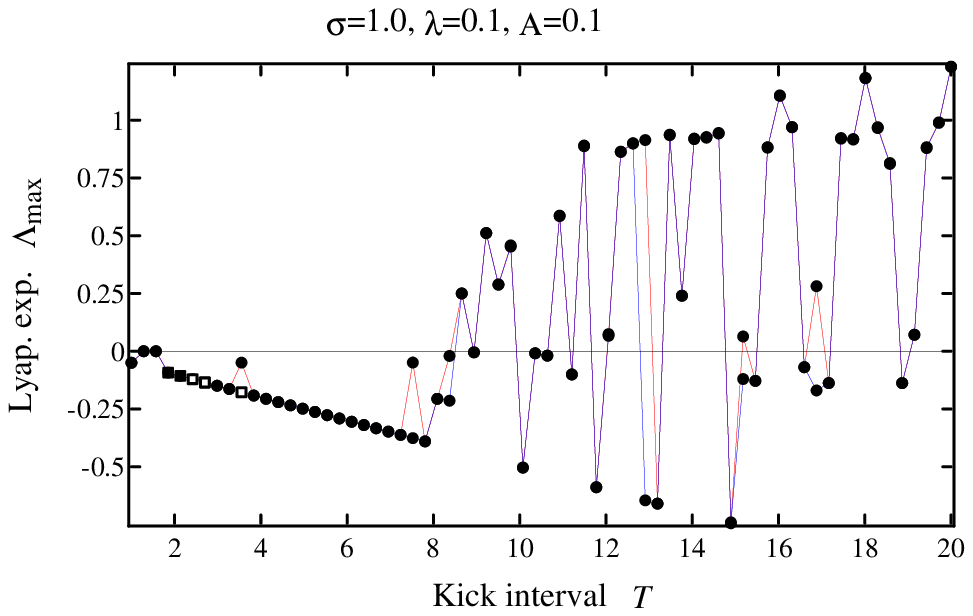}\vspace{18pt}\\
\includegraphics[bb=0in 0in 4in 2.5in]{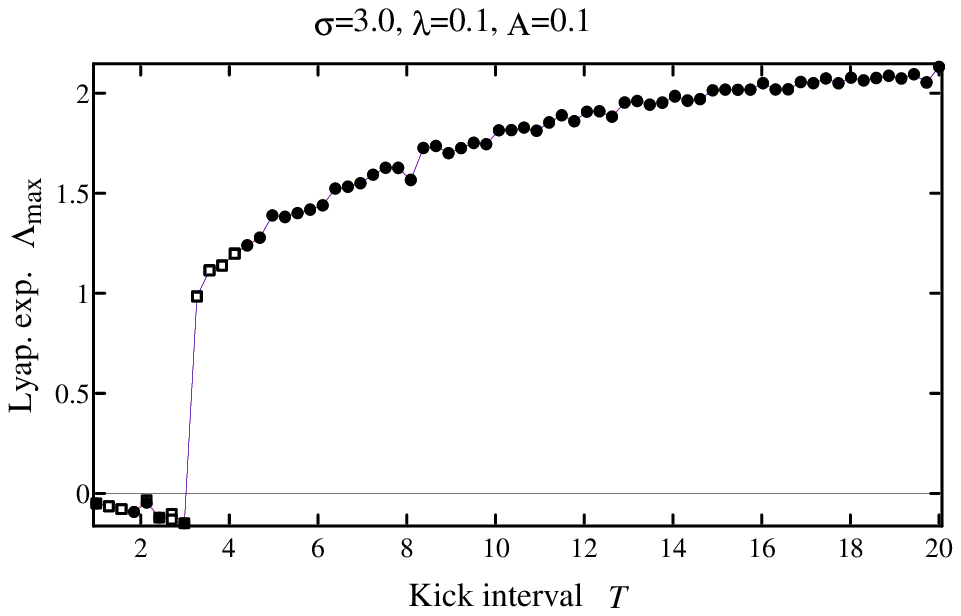}\\
\end{center}
\kaption{Effect of increasing shear on the Lyapunov exponents of
  the periodically-kicked linear shear flow.  Squares indicate
  that the corresponding $F_T$-orbit has veered outside the
  region $|y| < 0.15$.  Upper and lower estimates of $\lmax$ are
  both shown (see Simulation Details).}
\label{fig:kick-cycle}
\end{figure}

\begin{figure}
\begin{center}
\includegraphics[bb=0in 0in 4in 2.5in]{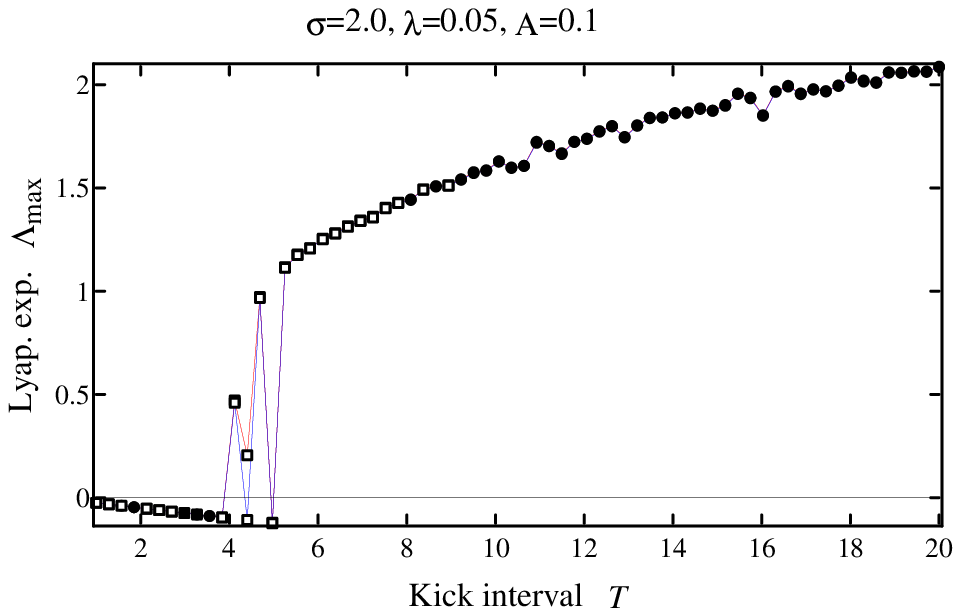}\vspace{18pt}\\
\includegraphics[bb=0in 0in 4in 2.5in]{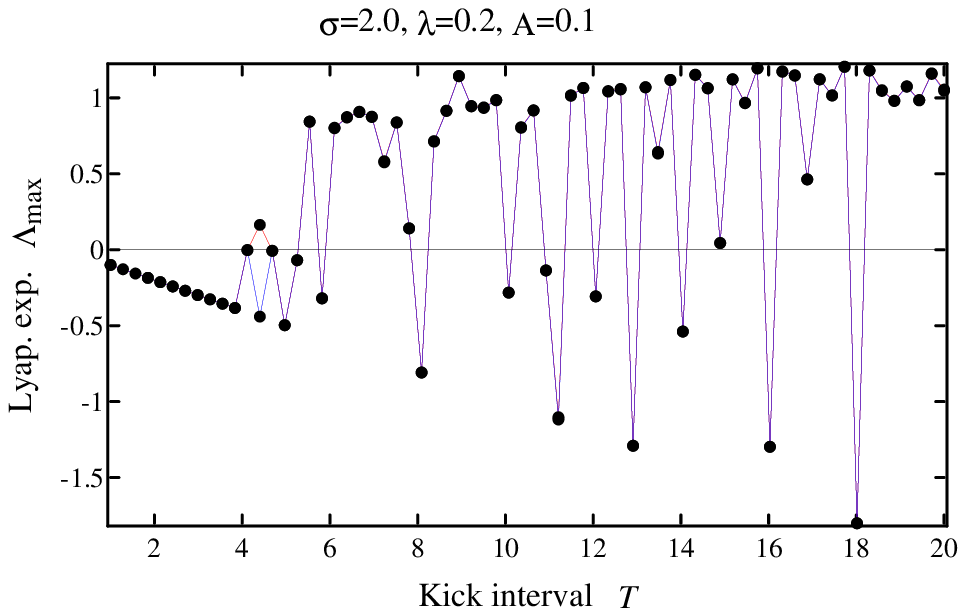}\vspace{18pt}\\
\includegraphics[bb=0in 0in 4in 2.5in]{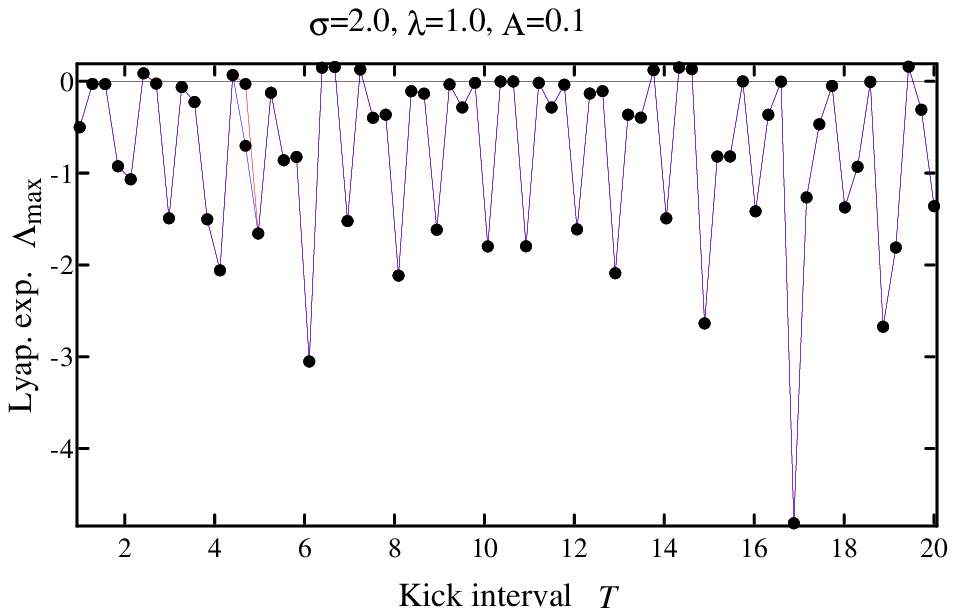}\\
\end{center}
\kaption{Effect of increasing damping on the Lyapunov exponents
  of the periodically-kicked linear shear flow.  Squares
  indicate that the corresponding $F_T$-orbit has veered outside
  the region $|y| < 0.15$.}
\label{fig:kick-cycle-too}
\end{figure}

\medskip
\noindent
{\em Supporting Numerical Evidence.} Figures
~\ref{fig:kick-cycle} and \ref{fig:kick-cycle-too} show the
largest Lyapunov exponent $\Lambda_{\max}$ of $F_T$ versus the
kick period $T$.  (Note that this is the expansion rate per kick
period and is $T$ times the rate per unit time.)  In
Fig.~\ref{fig:kick-cycle}, $\lambda$ and $A$ are fixed, and
$\sigma$ is increased.  We purposefully start with too small a
$\sigma$ so that we may see clearly the gradual changes in
$\Lambda_{\max}$.  The results are in excellent agreement with
the description at the end of Sect.~\ref{sec1.2} (which pertains
to regimes with very large $T$), even though $T$ is not so large
here: In the top picture, where $\frac{\sigma}{\lambda} A$ is
small, the plot confirms a competition between quasi-periodicity
and sinks; in the middle picture, we see first $\Lambda_{\max}$
becoming increasingly negative, then transitions into a
competition between transient and sustained chaos, with the
latter dominating in the bottom picture.
Fig.~\ref{fig:kick-cycle-too} shows the same phenomena in
reverse order, with $\sigma$ and $A$ fixed and $\lambda$
increasing.  Notice that even for $\sigma, \lambda$ and $A$
leading to chaotic dynamics, $\Lambda_{\max}$ is negative for
small $T$.  This is in agreement with the influence of the
factor $(1 - e^{-\lambda T})$ in Eq.~(\ref{time-t map}).

As explained in (a) above, when $\lambda T$ is too small
relative to $A$, orbits stray farther from $\gamma$.  Data
points corresponding to parameters for which this happens are
marked by open squares.  For purposes of demonstrating the
phenomena in question, there is nothing wrong with these data
points, but as explained earlier, caution must be exercised with
these parameters in systems where the basin of $\gamma$ is
smaller.

\bigskip
\noindent
{\em Simulation Details.} The numbers $\Lambda_{\max}$ are
computed by iterating the map in Eq.~(\ref{time-t map}) and its
Jacobian, and tracking the rate of growth of a tangent vector.
We use $4\times 10^5$ iterates of $F_T$ in each run.  Mindful of
the delicate situation due to competition between transient and
sustained chaos, and to lower the possibility of atypical
initial conditions, we perform 10 runs for each choice of
$(\sigma,A,\lambda,T)$, using for each run an independent,
random (with uniform distribution) initial condition
$(\theta_0,y_0)\in [0,1)\times[-0.1,0.1]$.  Among the 10 values
  of $\lmax$ computed, we discard the largest and the smallest,
  and plot the maximum and minimum of the remaining 8.  As one
  can see in Figs.~\ref{fig:kick-cycle} and
  \ref{fig:kick-cycle-too}, the two estimates occasionally do
  not agree.  This may be because not all initial conditions in
  the system have identical Lyapunov exponents, or it may be
  that the convergence to the true value of $\lmax$ is
  sufficiently slow and more iterates are needed, {\em i.e.}
  there are long transients.

%%%%%%%%%%%%%%%%%%%%%%%%%%%%%%%%%%%%%
%%%%%%%%%%%%%%%%%%%%%%%%%%%%%%%%%%%%%
\section{Study 2: Poisson Kicks}
\label{study2}

We consider next a variant of Eq.~(\ref{model1}) in which
deterministic, periodic kicks are replaced by ``random kicks.''
Here, random kicks refer to kicks at random times and with
random amplitudes.  More precisely, we consider
\begin{align}
\label{poisson}
\dot{\theta} &= 1 + \sigma y\\
\dot{y} &= -\lambda y + \sin(2\pi\theta)\sum_n \bA_n\delta(t-\bT_n)\nonumber
\end{align}
where the kick times $\bT_n$ are such that $\bT_{n+1}-\bT_n$,
$n=0,1,2,\cdots,$ are independent exponential random
variables with mean $T$, and the kick
amplitudes $\bA_n$ are independent and uniformly distributed
over the interval $[0.8\ A, 1.2\ A]$ for some $A > 0$.  (We do
not believe detailed properties of the laws of $\bT$ and $\bA$
have a significant impact on the phenomena being addressed.)
The analog here of the time-$T$ map in Study 1 is the {\it
  random map} $F = \Phi_\bT\circ K_\bA$ where $\bT$ and $\bA$
are random variables.

By the standard theory of random maps, Lyapunov exponents with
respect to stationary measures are well defined and are
nonrandom, {\it i.e.} they do not depend on the sample path
taken {\cite{kifer}}.  Notice that if $\sigma \neq 0$, the
system (\ref{poisson}) has a unique stationary measure which is
absolutely continuous with respect to Lebesgue measure on
$S^1\times{\mathbb R}$: starting from almost every $z_0 \in S^1
\times {\mathbb R}$, after one kick, the distribution acquires a
density in the $y$-direction; since vertical lines become
slanted under $\Phi_t$ due to $\sigma \neq 0$, after a second
kick the distribution acquires a (two-dimensional) density.

In terms of overall trends, our
 assessment of the likelihood of chaotic behavior 
 follows the analysis in Study 1 and will not be
repeated.  We identify the following two important differences:
\begin{itemize}

\item[(a)] {\em Smooth dependence on parameters.} Due to the
  averaging effects of randomness, we expect Lyapunov exponents
  to vary smoothly with parameter, without the wild oscillations
  in the deterministic case.

\item[(b)] {\em Effects of large deviations.}  A large number
  of kicks occurring in quick succession may have the following
  effects:
  \begin{itemize}

  \item[(i)] They can cause some orbits to stray far away from
    $\gamma = S^1\times\{0\}$.  This is guaranteed to happen, though
    infrequently, in the long run.  Thus, it is reasonable to
    require only that a large fraction --- not all --- of the
    stationary measure (or perhaps of the {\it random
      attractors} $\Gamma_\omega$) to lie in a prescribed
    neighborhood of $\gamma$.

  \item[(ii)] It appears possible, in principle, for a rapid
    burst of kicks to lead to chaotic behavior even in
    situations where the shear is mild and kick amplitudes are
    small.  To picture this, imagine a sequence of kicks sending
    (or maintaining) a segment far from $\gamma$, allowing the
    shear to act on it for an uncharacteristically long time.
    One can also think of such bursts as effectively setting
    $\lambda$ to near $0$ temporarily, creating a very
    large $\frac{\sigma}{\lambda}A$.  On the other hand, if
    $\sigma$ is small, then other forces in the system may try
    to coax the system to form sinks between these infrequent
    events.  We do not have the means to assess which scenario
    will prevail.

  \end{itemize}
\end{itemize}

\paragraph{Summary of Findings.} {\it In terms of overall trends,
the results are consistent with those in Study 1. Two
differences are observed. One is the rapid convergence of
$\Lambda_{\max}$ and their smooth dependence on parameters. The
other is that positive Lyapunov exponents for $F$ are found both
for smaller values of $\frac{\sigma}{\lambda} A$ and for
apparently very small $T$ (which is impossible for periodic
kicks), lending credence to the scenario described in (b)(ii)
above.}

\begin{figure}
\begin{center}
\includegraphics[bb=0in 0in 2.2in 2.2in]{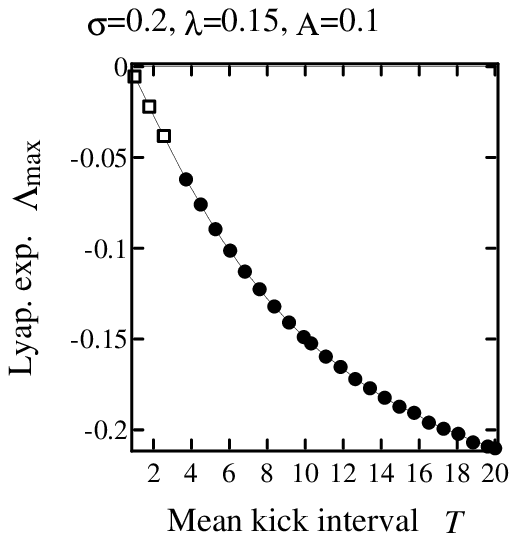}%
\includegraphics[bb=0in 0in 2.2in 2.2in]{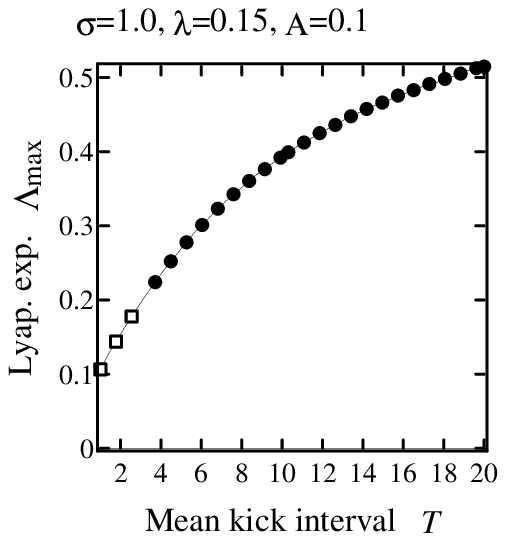}%
\includegraphics[bb=0in 0in 2.2in 2.2in]{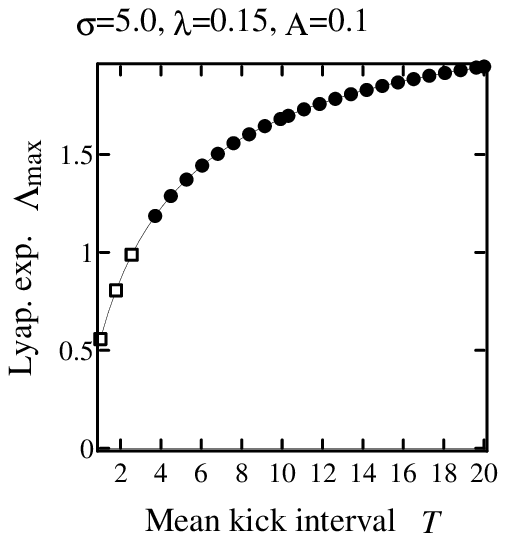}\\
(a) {\em Increasing shear}\\\vspace{18pt}
\includegraphics[bb=0in 0in 2.2in 2.2in]{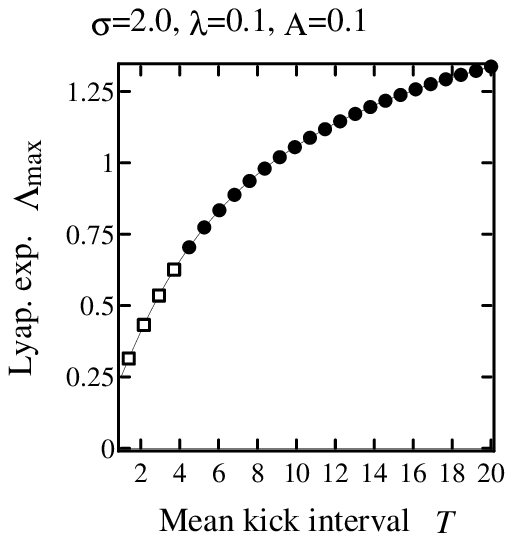}%
\includegraphics[bb=0in 0in 2.2in 2.2in]{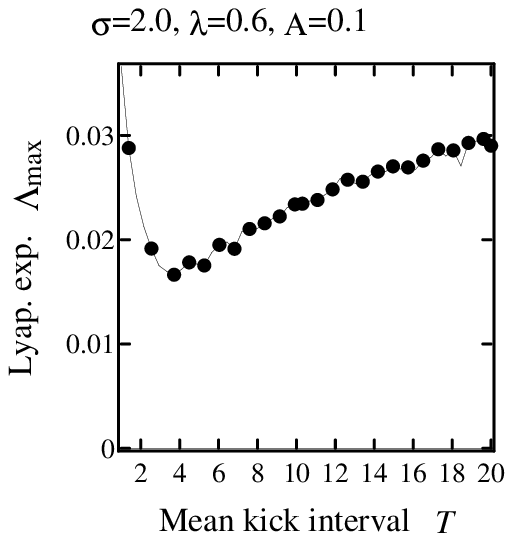}%
\includegraphics[bb=0in 0in 2.2in 2.2in]{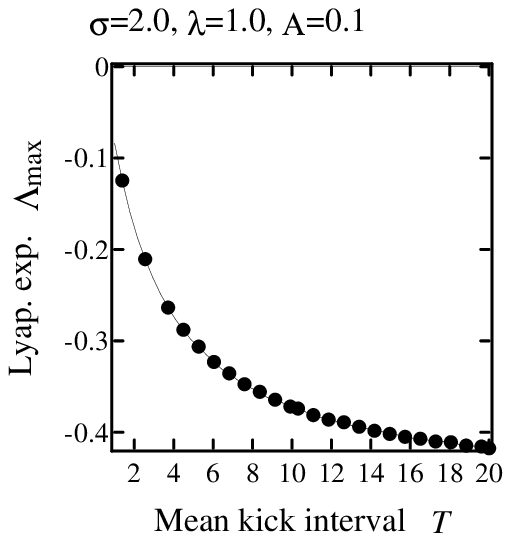}\\
(b) {\em Increasing damping}\\\vspace{18pt}
\end{center}
\kaption{Lyapunov exponents for the linear shear flow with
  Poisson kicks.  Squares indicate the corresponding orbit
  spends more than 20\% of the time in the region $|y| > 0.1$.}
\label{fig:poisson}
\end{figure}

\bigskip
\noindent {\em Supporting Numerical Evidence.}
Fig.~\ref{fig:poisson} shows $\Lambda_{\max}$ as a function of
the mean kick interval $T$.  As in Study 1, we first show the
effects of increasing $\sigma$ and then the effects of
increasing $\lambda$. Without the oscillations seen previously,
the present plots are straightforward to interpret. In case one
wonders how $\lmax$ curves can switch from strictly-decreasing
to strictly-increasing behavior, the middle panel of
Fig.~\ref{fig:poisson}(b) catches such a swtich ``in the act.''
Squares indicate that the orbit computed spends $>20\%$ of its
time outside of the region $\{|y|<0.1\}$.

%%%%%%%%%%%%%%%%%%%%%%%%%%%%%%%%%%%%%
%%%%%%%%%%%%%%%%%%%%%%%%%%%%%%%%%%%%%

\section{Study 3: Continuous-Time Stochastic Forcing}
\label{study3}

In this section, we investigate the effect of forcing by white
noise.  The resulting systems are described by stochastic differential
equations (SDEs).  We consider two ways to force the system:

\bigskip
\noindent {\bf Study 3a: Degenerate white noise applied in
chosen direction:}
\begin{align}
\label{deg}
d\theta &= (1 + \sigma y)\ dt\\
dy &= -\lambda y\ dt + a\sin(2\pi\theta) \ dB_t\nonumber
\end{align}
\smallskip
\noindent {\bf Study 3b: Isotropic white noise:}
\begin{align}
\label{nondeg}
d\theta &= (1 + \sigma y)\ dt  + a \sin(2\pi\theta)\ dB^1_t\\
dy &= -\lambda y\ dt + a\sin(2\pi\theta)\  dB^2_t\nonumber
\end{align}
In Study 3a, $B_t$ is standard 1-dimensional Brownian motion
(meaning with variance $=1$). In Study 3b, $(B^1_t, B^2_t)$ is a
standard 2-D Brownian motion, {\it i.e.}, they are independent
standard 1-D Brownian motions.  For definiteness, we assume the
stochastic terms are of It\^{o} type.  Notice that the two
parameters $A$ and $T$ in Studies 1 and 2 have been combined
into one, namely $a$, the coefficient of the Brownian noise.

\bigskip

By standard theory {\cite{arnbook,kunita}}, the solution process
of an SDE can be represented as a stochastic flow of
diffemorphisms. More precisely, if the coefficients of the SDE
are time-independent, then for any time step $\Delta t>0$, the
solution may be realized, sample path by sample path, as the
composition of random diffeomorphisms \ $\cdots \ \circ f_3
\circ f_2 \circ f_1$, where the $f_i$ are chosen {\it
  i.i.d.}\ with a law determined by the system (the $f_i$ are
time-$\Delta t$ flow-maps following this sample path).  This
representation enables us to treat an SDE as a {\it random
  dynamical system} and to use its Lyapunov exponents as an
indicator of chaotic behavior.  It is clear that system
(\ref{nondeg}) has a unique invariant density, which is the
solution of the Fokker-Planck equation.  Even though the
stochastic term in system (\ref{deg}) is degenerate, for the
same reasons discussed in Study 3, it too has a unique
stationary measure, and this measure has a density.  The
Lyapunov exponents considered in this section are with respect
to these stationary measures.

Before proceeding to an investigation of the two systems above, we first
comment on the case of purely additive noise, {\it i.e.}
Eq.~(\ref{nondeg}) without the $\sin(2\pi\theta)$ factor in
either Brownian term.  In this case it is easy to see that 
all Lyapunov exponents are $\leq 0$, for the random maps are
approximately time-$\Delta t$ maps of the unforced flow composed
with random (rigid) translations.  Such a system is clearly not chaotic.

With regard to system (\ref{deg}), we believe that even
though the quantitative estimates from Study 1 no longer apply, 
a good part of the {\it qualitative reasoning} behind the arguments 
continues to be valid.  In
particular, we conjecture that
\begin{itemize}

\item[(a)] trends, including qualitative dependences on $\sigma$
  and $\lambda$, are as in the previous two studies;

\item[(b)] the effects of large deviations noted for Poisson
  kicks (Study 2, item (b)) are even more prominent here, given
  that the forcing now occurs continuously in time.

\end{itemize}

As for system (\ref{nondeg}), we expect it to be less effective
in producing chaos, {\it i.e.} more inclined to form sinks, than
system (\ref{deg}).  This expectation is based on the following
reasoning: Suppose first that we force {\it only} in the
$\theta$-direction, {\it i.e.}, suppose the $dB^2_t$ term in
(\ref{nondeg}) is absent.  Then the stochastic flow leaves
invariant the circle $S^1\times\{0\}$, which is the limit cycle
of the deterministic part of the system.  A general theorem
tells us that when a random dynamical system on a circle has an
invariant density, its Lyapunov exponent is always $\leq 0$; in
this case, it is in fact strictly negative because of the
inhomogeneity caused by the sine function {\cite{kifer}}.  Thus
the corresponding 2-D system has ``random sinks.''  Now let us
put the $y$-component of the forcing back into the system.  We
have seen from previous studies that forcing the $y$ direction
alone may lead to chaotic behavior.  The tendency to form sinks
due to forcing in the $\theta$-direction persists, however, and
weakens the effect of the shear-induced stretching.

We now discuss the results of simulations performed to validate
these ideas.

%%%%%%%%%%%%%%%%%%%%%%%%%%%%%%%%
%\subsection{Results of simulations}

\paragraph{Summary of Findings.} {\it 
\begin{itemize}

\item[(i)] In the case of {\rm degenerate} white noise, the
  qualitative dependence of $\Lambda_{\max}$ on $\sigma$ and
  $\lambda$ are as expected, and the effects of large deviations
  are evident. In particular, $\Lambda_{\max}$ is positive for
  very small values of $\sigma,\lambda$ and $a$ provided
  $\frac{\sigma}{\lambda}$ is large. {\rm This cannot happen for
    periodic kicks; we attribute it to the effect of large
    deviations.}

\item[(ii)] Isotropic white noise
is considerably less effective in producing chaos than 
forcing in the $y$-direction only, meaning
it produces a smaller (or more negative) $\Lambda_{\max}$.

\item[(iii)] In both cases, we discover the following approximate
  scaling: Under the scaling transformations $\lambda\mapsto
  k\lambda, \sigma\mapsto k\sigma$ and $a\mapsto \sqrt{k}a$,
  $\Lambda_{\max}$ transforms approximately as
  $\Lambda_{\max}\mapsto k\Lambda_{\max}$. In the case of
  degenerate white noise, when both $\sigma$ and $\frac{\lambda}{\sigma}$ 
  are not too small ({\it e.g.}, $> 3$), this scaling gives excellent
  predictions of $\Lambda_{\max}$ for the values computed.
  \end{itemize}}

\medskip We remark that (iii) does not follow by scaling time in
the SDE.  Indeed, scaling time by $k$ in Eq.~(\ref{deg}), we
obtain
\begin{align}
\label{scaleddeg}
d\theta &= (k + k \sigma y)\ dt \ , \\
dy &= - k \lambda y\ dt + \sqrt k a\sin(2\pi\theta)\  dB_t\ .\nonumber
\end{align}
Thus the approximate scaling in (iii) asserts that the Lyapunov
exponent of system (\ref{scaleddeg}), equivalently $k$ times the
$\Lambda_{\max}$ for Eq.~(\ref{deg}), is roughly equal to that
of the system obtained by changing the first equation in
(\ref{scaleddeg}) to $d\theta = (1+ k \sigma y)dt$. In other
words, $\Lambda_{\max}$ seems only to depend minimally on the
frequency of the limit cycle in the unforced system.

\begin{figure}[t]
\begin{center}
\includegraphics[bb=0in 0in 2.25in 2.25in]{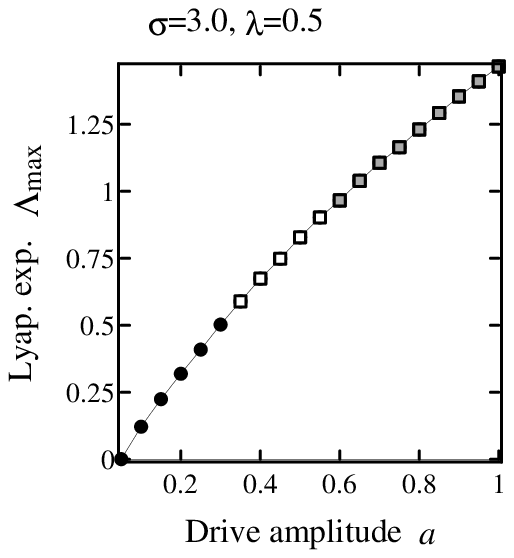}%
\includegraphics[bb=0in 0in 2.25in 2.25in]{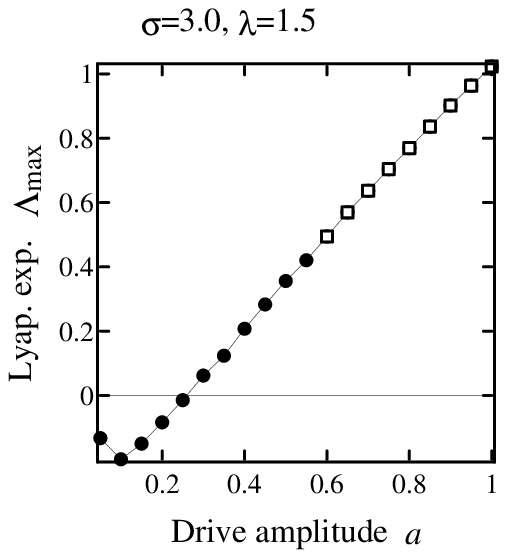}%
\includegraphics[bb=0in 0in 2.25in 2.25in]{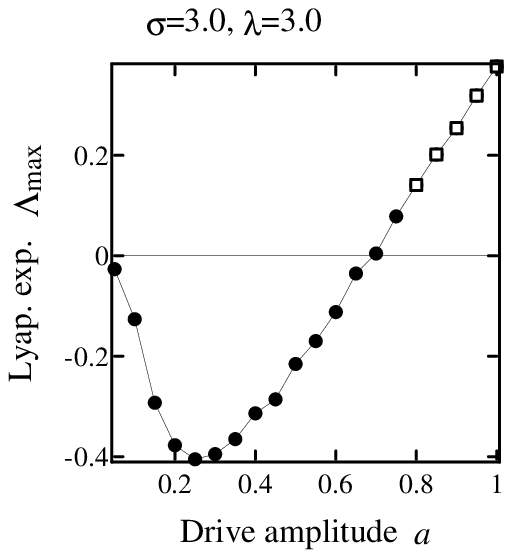}\\
\end{center}
\kaption{Lyapunov exponents for the linear shear flow driven by
  degenerate white noise (Eq.~(\ref{deg})).  Open squares
  indicate that the corresponding orbits spend more than 20\% of
  the time in the region $|y| > 0.3$; shaded squares do the same
  for the region $|y| > 0.5$.}
\label{fig:whitenoise}
\end{figure}

\begin{figure}[h]
\begin{center}
\includegraphics[bb=0in 0in 2.75in 2.25in]{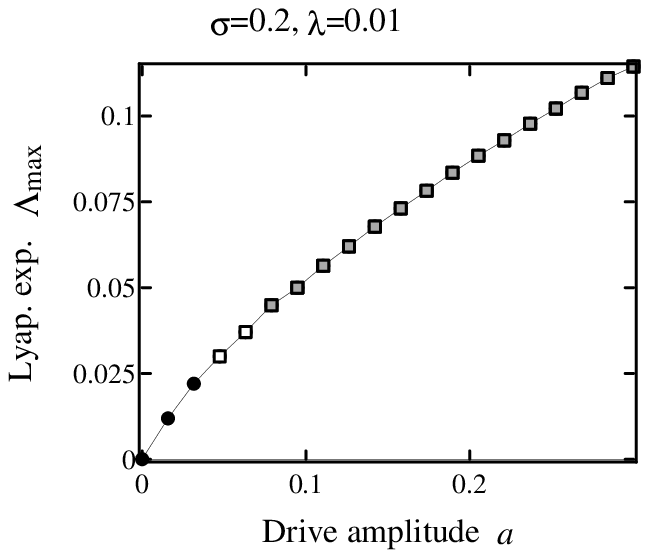}
\end{center}
\kaption{Lyapunov exponents for the linear shear flow driven by
  degnerate white noise, for small values of $\sigma$ and
  $\lambda$.}
\label{fig:small sigma}
\end{figure}

\begin{figure}[h]
\begin{center}
\vspace{12pt}\includegraphics[bb=0in 0in 2.25in 2.25in]{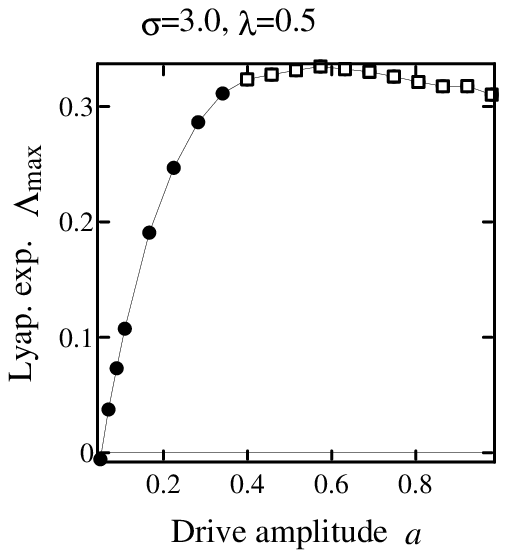}%
\includegraphics[bb=0in 0in 2.25in 2.25in]{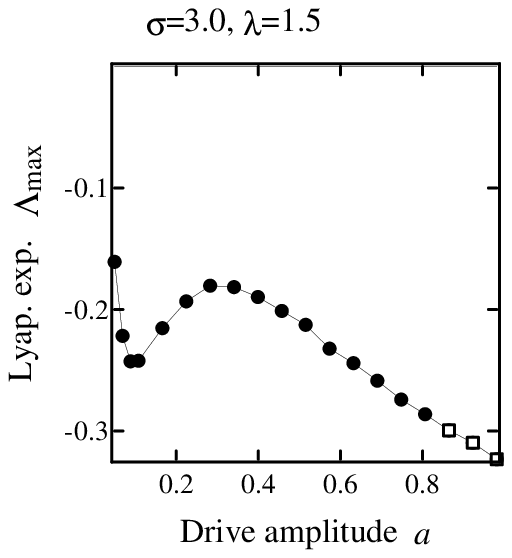}%
\includegraphics[bb=0in 0in 2.25in 2.25in]{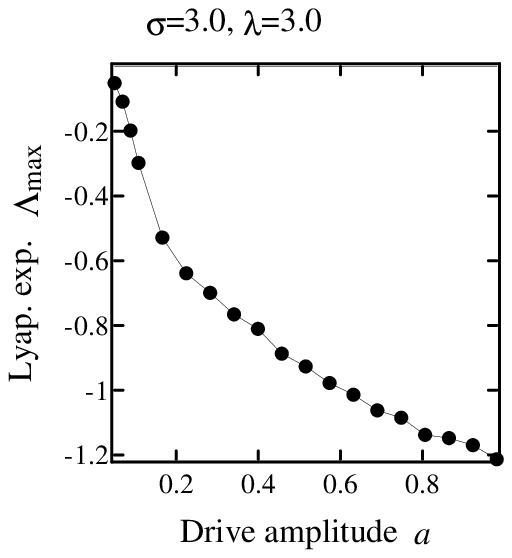}\\
(a) $\sigma = 3.0$\\
\vspace{12pt}\includegraphics[bb=0in 0in 2.25in 2.25in]{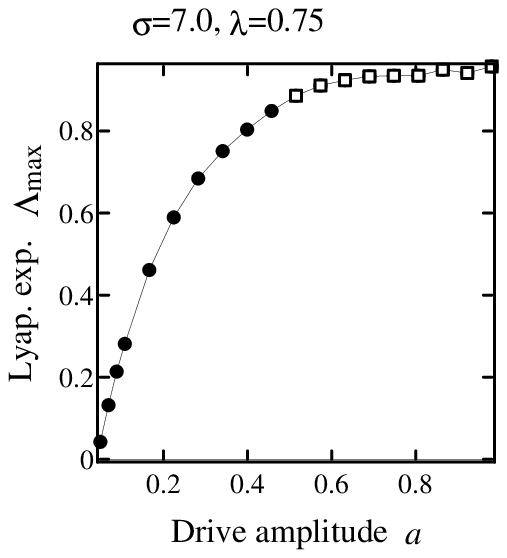}%
\includegraphics[bb=0in 0in 2.25in 2.25in]{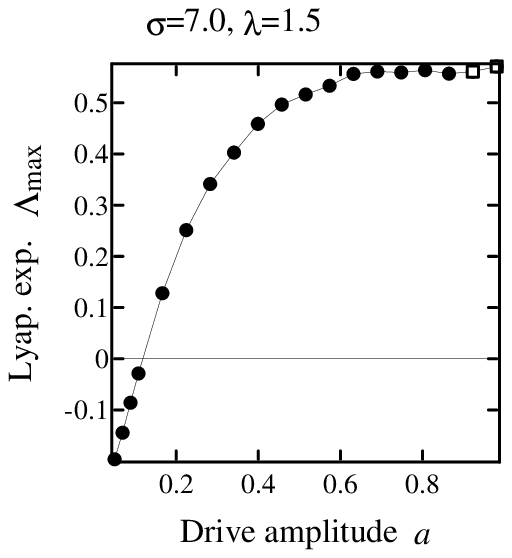}%
\includegraphics[bb=0in 0in 2.25in 2.25in]{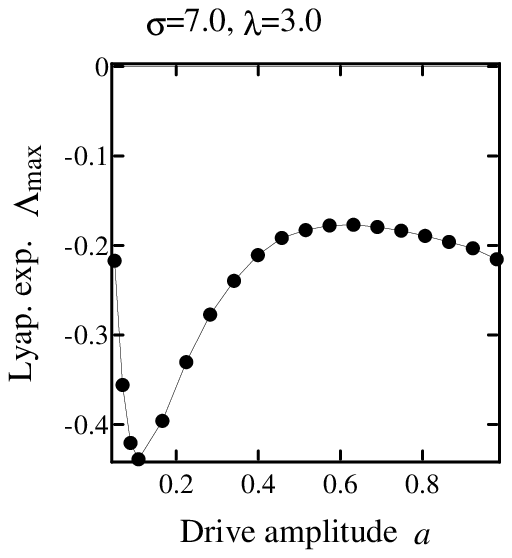}\\
(b) $\sigma = 7.0$
\end{center}
\kaption{Lyapunov exponent for the linear shear flow driven by
  isotropic white noise (Eq.~(\ref{nondeg})).  Squares indicate
  that the corresponding orbits spend more than 20\% of the time
  in the region $|y| > 0.3$.}

\label{fig:isotropic}
\end{figure}

\begin{figure}
\begin{center}
\includegraphics[bb=0in 0in 3.5in 2.5in]{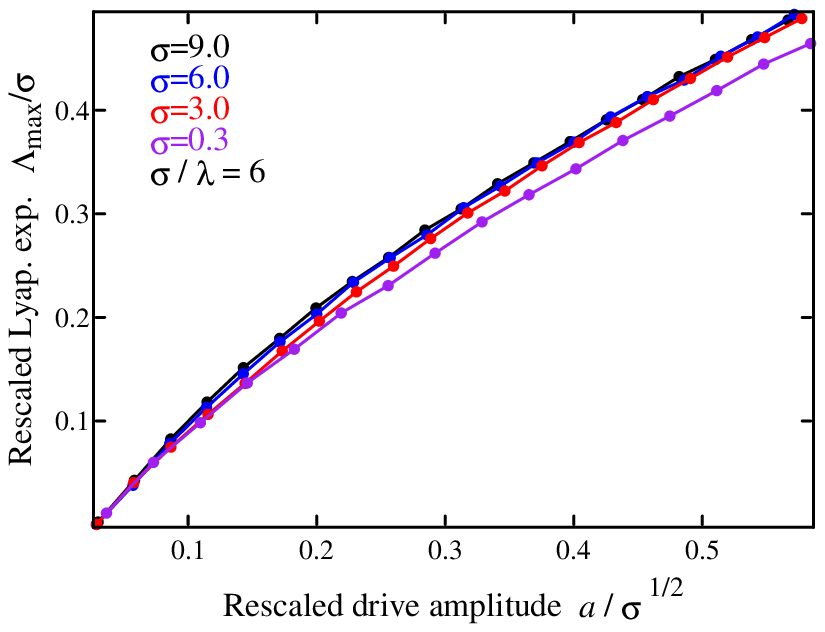}
\end{center}
\kaption{Evidence of scaling: We fix $\frac{\sigma}{\lambda}$
  and plot $\lmax/\sigma$ as functions of the rescaled drive
  amplitude $a/\sqrt{\sigma}$; from top to bottom, the curves
  are in order of decreasing $\sigma$.}
\label{fig:scaling}
\end{figure}

\bigskip
\noindent {\em Supporting Numerical Evidence.}  Plots of
$\Lambda_{\max}$ as functions of $a$ are shown in
Figs.~\ref{fig:whitenoise} -- \ref{fig:scaling}. 

In Fig.~\ref{fig:whitenoise}, the forcing is degenerate, and for
fixed $\sigma$, $\lmax$ decreases with increasing damping as
expected.  Notice that compared to the two previous studies, a
somewhat larger damping is required to maintain a good fraction
of the attractor near $\gamma$.

Fig.~\ref{fig:small sigma} shows that $\Lambda_{\max}$ is
positive for values of $\sigma$ and $\lambda$ as small as $0.2$
and $0.01$, and white noise amplitudes $a$ close to $0$.  Notice
first that this is consistent with the scaling conjectured in
(iii) above, and second that in the case of periodic kicks,
comparable values of $\sigma$ and $\lambda$ would require a
fairly substantial kick, not to mention long relaxation periods,
before chaotic behavior can be produced.  We regard this as
convincing evidence of the significant effects of large
deviations in continuous-time forcing.  (It must be pointed out,
however, that in our system, the basin of attraction is the
entire phase space, and a great deal of stretching is created
when $|y|$ is large. That means system (\ref{deg}) takes greater
advantage of large deviations than can be expected ordinarily.

Fig.~\ref{fig:isotropic}(a) shows $\Lambda_{\max}$ in the isotropic case 
for the same parameters as in Fig.~\ref{fig:whitenoise}. A comparison
of the two sets of results confirms the conjectured tendency toward 
negative exponents when the forcing is isotropic.  
Fig.~\ref{fig:isotropic}(b) shows that this
tendency can be overcome by increasing $\sigma$.

Fig.~\ref{fig:scaling} shows four sets of results, overlaid on
one another, demonstrating the scaling discussed in item (iii)
above. Fixing $\frac{\sigma}{\lambda}=6$, we show the graphs of
$\Lambda_{\max}/\sigma$ as functions of $a/\sqrt{\sigma}$ for
four values of $\sigma$. The top two curves (corresponding to
$\sigma=6$ and $9$) coincide nearly perfectly. Similar
approximate scalings, less exact, are observed for smaller
values of $\frac{\sigma}{\lambda}$, both when $\Lambda_{\max}$
is positive and negative.

\bigskip
\noindent {\em Simulation Details.}  We compute Lyapunov
exponents numerically by solving the corresponding variational
equations (using an Euler solver with time steps of $10^{-5}$)
and tracking the growth rate of a tangent vector.  To account
for the impact of the realization of the forcing on the computed
exponents, for each choice of $(\sigma,\lambda,a)$ we perform 12
runs in total, using 3 independent realizations of the forcing
and, for each realization, 4 independent initial conditions
(again uniformly-distributed in $[0,1)\times[-0.1,0.1]$).  For
  almost all the parameter values, the estimates agree to fairly
  high accuracy, so we simply average over initial conditions
  and plot the result.

\paragraph{Related Results.}
The asymptotic stability of dynamical systems driven by random
forcing has been investigated by many authors using both
numerical and analytic methods. Particularly relevant to our
study are results pertaining to the random forcing of
oscillators (such as Duffing-van der Pol oscillators) and
stochastic Hopf bifurcations; see e.g.
{\cite{arn,aus,bax02,bax03,bax04,bax-gouk,pinsky,nsri}}.  Most
of the existing results are perturbative, {\it i.e.}, they treat
regimes in which both the noise and the damping are very small.
Positive Lyapunov exponents are found under certain conditions.
We do not know at this point if the geometric ideas of this
paper provide explanations for these results.

%%%%%%%%%%%%%%%%%%%%%%%%%%%%%%%%%%
%%%%%%%%%%%%%%%%%%%%%%%%%%%%%%%%%%%
\section{Study 4: Sheared-Induced Chaos in Quasiperiodic Flows}
\label{study4}

\heading{Model and Background Information}

In this section, we will show that external forcing can lead to
shear-induced chaos in a coupled phase oscillator system of the
form
\begin{center}
\input{pix/oscillators.pic}
\end{center}
\noindent The governing equations are
\begin{align}
\label{eq:2osc}
\dot{\theta}_1 &= \nu_1 + z(\theta_1) [\aff g(\theta_2) + I(t)]\ ,\\
\dot{\theta}_2 &= \nu_2 + z(\theta_2) [\afb g(\theta_1)].\nonumber
\end{align}
The state of the system is specified by two angles,
$(\theta_1,\theta_2)$, so that the phase space is the torus
${\mathbb T}^2\equiv [0,1)^2$.  The constants $\nu_1$ and
  $\nu_2$ are the oscillators' intrinsic frequencies; we set
  $\nu_1=1$ and $\nu_2=1.1$ (representing similar but not
  identical frequencies).  The constants $\aff$ and $\afb$
  govern the strengths of the feedforward and feedback
  couplings.  The oscillators are pulse-coupled: the coupling is
  mediated by a bump function $g$ supported on $[-\frac{1}{20},
    \frac{1}{20}]$ and normalized so that $\int_0^1
  g(\theta)\ d\theta = 1$.  The function $z(\theta)$, which we
  take to be $z(\theta)=\frac{1}{2\pi}[1-\cos (2\pi \theta)]$,
  specifies the sensitivity of the oscillators to perturbations
  when in phase $\theta$.  Finally, we drive the system with an
  external forcing $I(t)$, which is applied to only the first
  oscillator.  This simple model arises from neuroscience
  {\cite{winfree,taylor-holmes}} and is examined in more detail
  in {\cite{reliability}}.

Let $\Phi_t$ denote the flow of the unforced system, {\it i.e.},
with $I(t)\equiv 0$. Flowlines are roughly northeasterly and are
linear except in the strips $\{|\theta_1|<\frac{1}{20}\}$ and
$\{|\theta_2|<\frac{1}{20}\}$, where they are bent according to
the prescribed values of $\aff$ and $\afb$.  Let $\rho$ denote
the rotation number of the first return map of $\Phi_t$ to the
cross-section $\{\theta_2 = 0\}$.  It is shown in {\cite{prl}}
that for $\aff=1$, $\rho$ is monotonically increasing (constant
on extremely short intervals)
as one increases $\afb$, until it reaches $1$ at $\afb =
\afb^*\approx 1.4$, after which it remains constant on a large
interval.  At $\afb = \afb^*$, a limit cycle emerges in which
each oscillator completes one rotation per period; we say the
system is 1:1 phase-locked, or simply {\em phase-locked}.  
In {\cite{prl}}, it is shown
numerically that forcing the system by white noise after the
onset of phase-locking leads to $\Lambda_{\max} > 0$.  The
authors of {\cite{prl}} further cite Wang-Young theory
(the material reviewed in Sect.~\ref{wang-young}) as a
geometric explanation for this phenomenon.

In this section, we provide geometric and numerical evidence of
shear-induced chaos both before and after the onset of
phase-locking at $\afb=\afb^*$.  Our results for $\afb>\afb^*$
support the assertions in {\cite{prl}}. For $\afb<\afb^*$,
they will show that {\it limit cycles are not preconditions
  for shear-induced chaos}.  We will show that in
Eq.~(\ref{eq:2osc}), the mechanism for folding is already in
place before the onset of phase-locking, where the system is
quasi-periodic or has periodic orbits of very long periods; the
distinction between these two situations is immaterial since we
are concerned primarily with finite-time dynamics.  In the rest
of this section, we will, for simplicity, refer to the regime
prior to the onset of phase-locking as ``near-periodic.''

\heading{Folding: Geometric Evidence of Chaos}

The dynamical picture of kicks followed by a period of relaxation has a
simpler, more clear-cut geometry than that of continuous, random
forcing.  Thus we use the former to demonstrate why one may expect chaotic
behavior over the parameter ranges in question. The kick map is denoted
by $\kappa$ as in Section 1.

\bigskip
\noindent
{\em Folding in the periodic ({\em i.e.} phase-locked) regime.}
We will use $\afb=1.47$ for illustration purposes; similar
behavior is observed over a range of $\afb$ from $1.4$ to $1.6$.
Note that the system is phase-locked for a considerably larger
interval beyond $\afb=1.6$, but the strength of attraction grows
with increasing $\afb$, and when the attraction becomes too
strong, it is harder for folding to occur.

\begin{figure}
\begin{center}
\includegraphics[bb=0 0 269 221]{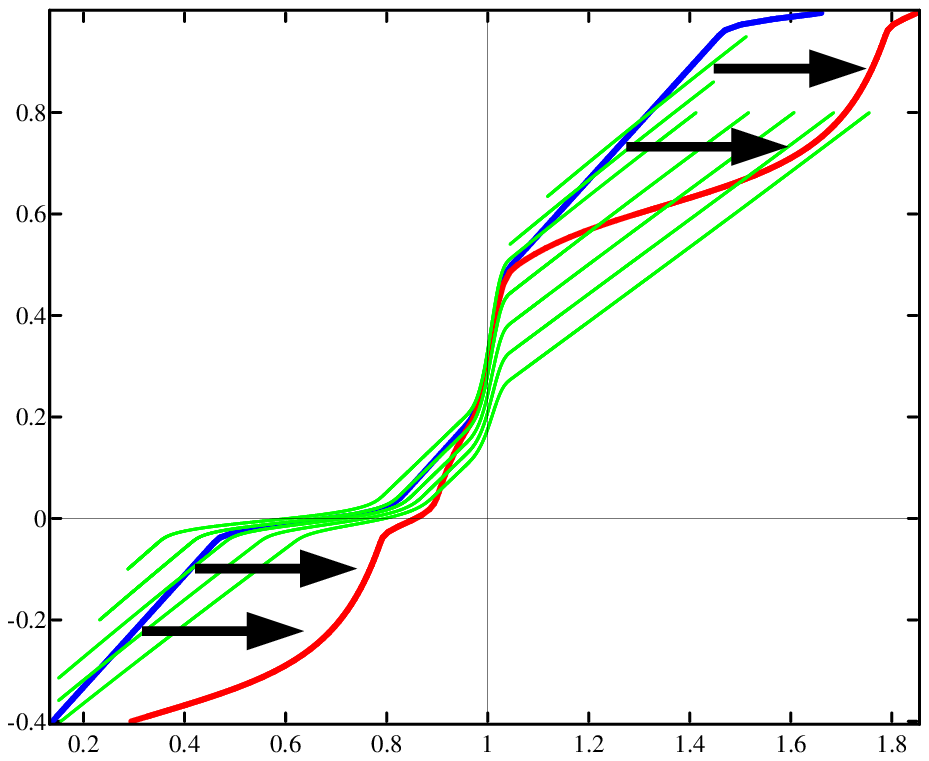}
\end{center}
\kaption{The strong-stable foliation of the system
  (\ref{eq:2osc}) in the phase-locked regime.  {\em Blue:} a
  lift of the limit cycle.  {\em Red:} the image of the cycle
  after a single kick.  {\em Green:} strong-stable foliation.
  Here, the parameters are $\nu_1=1$, $\nu_2=1.1$, $\aff=1$, and
  $\afb=1.47$.}
\label{cycle+wss}
\end{figure}

\begin{figure}[h]
\begin{center}
\includegraphics[bb=0 0 176 151]{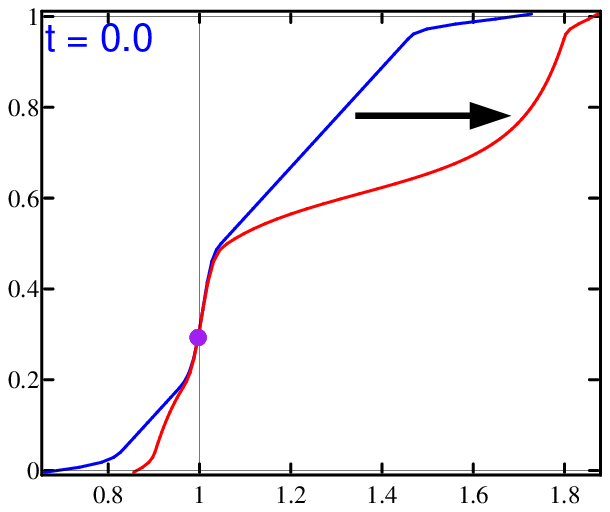}\vspace{18pt}\\
\includegraphics[bb=0in 0in 2.034277576757482in 2.5in]{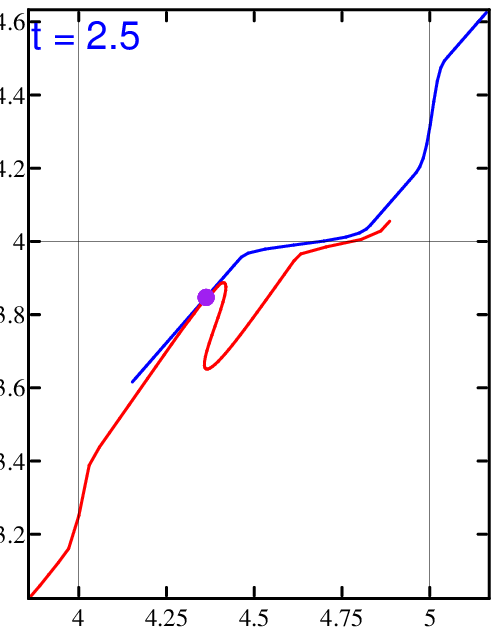}\hspace{24pt}%
\includegraphics[bb=0in 0in 2.2753509103534504in 2.5in]{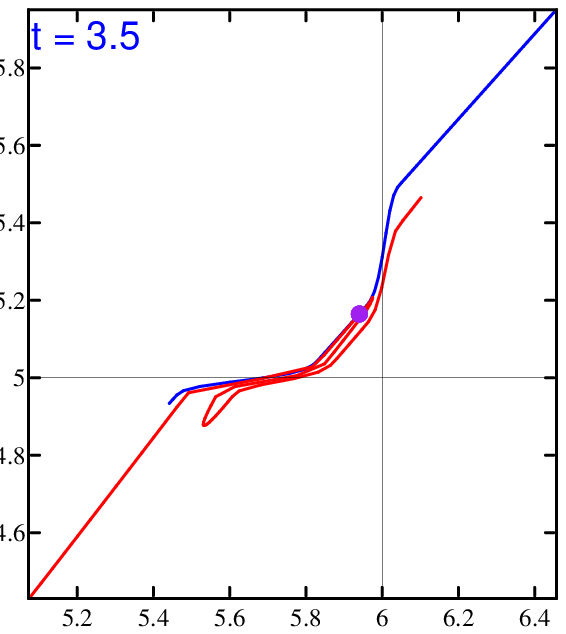}
\end{center}
\kaption{Snapshots of the limit cycle and its kicked image in a
  moving frame.  {\em Blue curves:} $\Phi_t$-images of
  $\gamma_0$, the part of the limit cycle between $\theta_2=0$
  and $\theta_2=1$.  {\em Red curves:} $\kappa(\gamma_0)$ and
  its images.  {\em Purple dot:} the point on $\gamma_0$ which
  does not move under $\kappa$.  The parameters are the same as
  in Fig.~\ref{cycle+wss}.}
\label{cycle snapshots}
\end{figure}

Fig.~\ref{cycle+wss} shows the limit cycle $\gamma$ ({\em blue
  curve}) of the unforced system at $\afb=1.47$; more precisely,
it shows a ``lift'' of $\gamma$ to ${\mathbb R}^2$, identifying
the torus ${\mathbb T}^2$ with ${\mathbb R}^2/{\mathbb Z}^2$.
Also shown is the image $\kappa(\gamma)$ of the cycle after a
single kick ({\em red curve}), where the kick map $\kappa$
corresponds to $I(t)= A \sum_n \delta(t-nT)$ with $A=1.5$, {\em
  i.e.}  $\kappa$ is given by $\kappa= \lim_{\eps \to 0}
\kappa_\eps(\eps)$ where $\kappa_\eps(t)$ is the solution of
$\dot{\theta}_1 = \frac{A}{\eps} z(\theta_1), \dot \theta_2=0$.
Notice the special form of the kicks: $\kappa$ acts
horizontally, and does not move points on $\theta_1=0$.  In
particular, $\kappa$ fixes a unique point $(0,b)$ on the cycle;
this point is, in fact, not affected by {\it any} kick of the
form considered in Eq. (\ref{eq:2osc}).  Several segments of
strong stable manifolds ({\em green curves}) of the unforced
system are drawn.  Recall that if $p$ is the period of cycle and
$n \in {\mathbb Z}^+$, then $\Phi_{np}(\kappa(z))$ lies on the
$W^{ss}$-curve through $\kappa(z)$ and is pulled toward the
cycle as $n$ increases (see Sect.~\ref{sec1.2}).  From the
relation between the $W^{ss}$-curves and the cycle, we see that
for $z\in\gamma$, $\Phi_t(\kappa(z))$ will lag behind
$\Phi_t(z)$ during the relaxation period.  Notice in particular
that there are points on $\kappa(\gamma)$ above the line
$\theta_2=b$ that are pulled toward the part of $\gamma$ below
$\theta_2=b$.  Since $(0,b)$ stays put, we deduce that some
degree of folding will occur if the time interval between kicks
is sufficiently long.

Fig.~\ref{cycle snapshots} illustrates how this folding happens
through three snapshots.  We begin with a segment
$\gamma_0\subset\gamma$ between $\theta_2=0$ and $\theta_2=1$
({\em blue curve}) and its image after a single kick ({\em red
  curve}). Both curves are then evolved forward in time and their images 
  at $t=2.5$ and $t=3.5$ are shown.  The purple dot marks the
point on $\gamma_0$ which does not move when kicked.  Notice
that these pictures are shown in a {\it moving frame} to
emphasize the geometry of $\Phi_t(\kappa(\gamma_0))$ relative to
$\Phi_t(\gamma_0)$.

\begin{figure}[h]
\begin{center}
\includegraphics[bb=0 4 156 154]{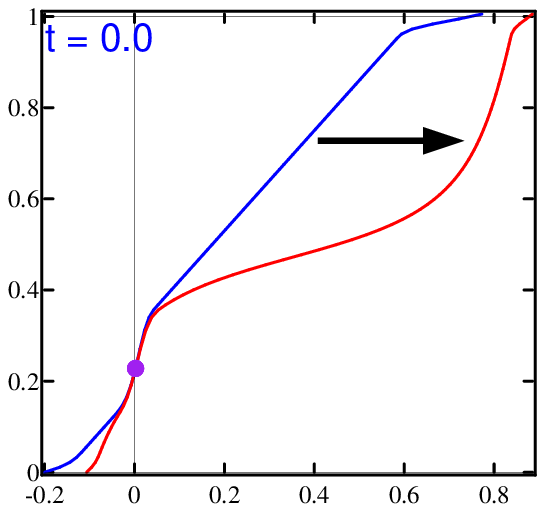}\hspace{24pt}%
\includegraphics[bb=0in 0in 2.1693296140636673in 2in]{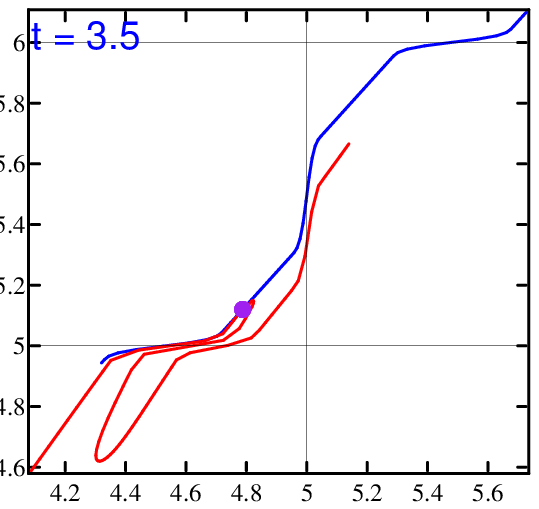}\\
\vspace{24pt}\includegraphics[bb=0in 0in 2.478107762552729in 2.5in]{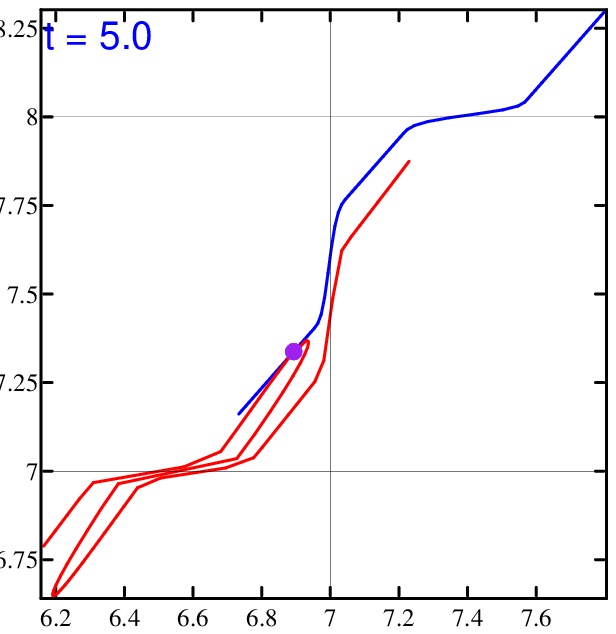}
\end{center}
\kaption{Snapshots of an orbit segment and its image after a
  single kick in a moving frame, for the system (\ref{eq:2osc})
  in a near-periodic regime.  {\em Blue curves:} a segment
  $\gamma_0$ of an orbit and its forward images
  $\Phi_t(\gamma_0)$ at $t=3.5$, $5$.  {\em Red curves:}
  $\kappa(\gamma_0)$ and its forward images.  {\em Purple dot:}
  the point on $\gamma_0$ which does not move under $\kappa$.
  The parameters are $\nu_1=1$, $\nu_2=1.1$, $\aff=1$, and
  $\afb=1.2$.}
\label{quasiper snapshots}
\end{figure}

\begin{figure}
\begin{center}
\includegraphics[bb=0in 0in 3.5in 3in]{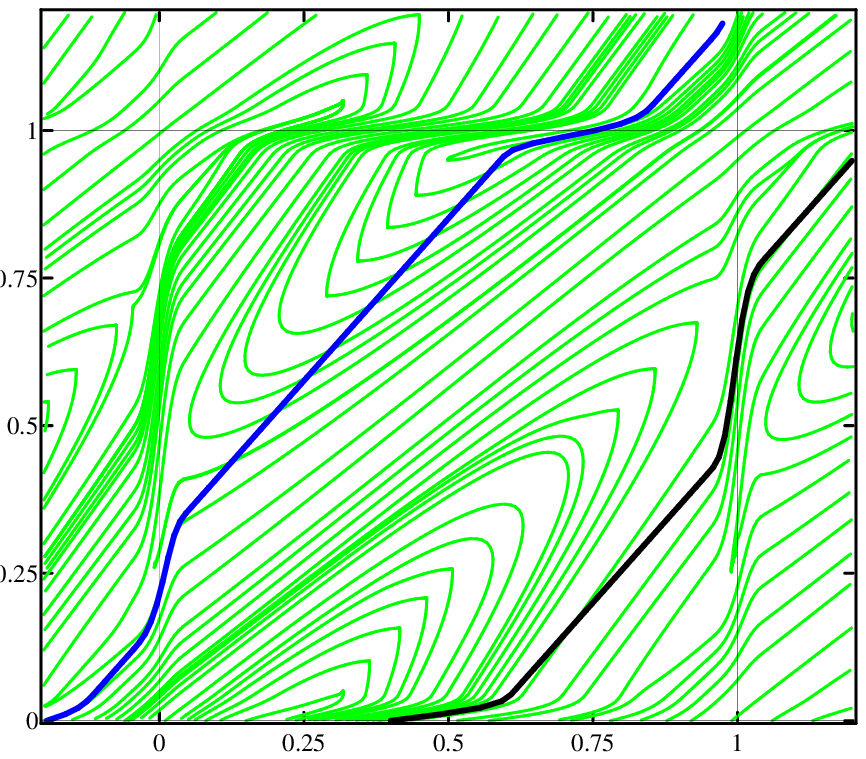}
\end{center}
\kaption{The time-5 stable foliation of the system
  (\ref{eq:2osc}) in a near-periodic regime.  {\em Blue \& black
    curves:} two orbit fragments.  {\em Green curves:} time-5
  stable foliation.  The parameters here are the same as in
  Fig.~\ref{quasiper snapshots}.}
\label{finite-time-wss}
\end{figure}

\bigskip
\noindent
{\em Folding in the near-periodic regime.}  Fig.~\ref{quasiper
  snapshots} shows snapshots of a similar kind for $\afb=1.2$;
this value of $\afb$ puts the system in the near-periodic
regime.  The snapshots begin with an (arbitrary) orbit segment
$\gamma_0$ and its image $\kappa(\gamma_0)$; the location of
$\gamma_0$ is near that of the limit cycle in
Fig.~\ref{cycle+wss}.  The kicked segment clearly folds; indeed,
the picture is qualitatively very similar to that of the limit
cycle case.  Note that at $\afb=1.2$, the rotation number of the
return map to $\{\theta_2=0\}$ is a little below 1, so that
$\Phi_t(\kappa(\gamma_0))$ has an overall, slow drift to the
left when viewed in the fixed frame $[0,1)^2$.  This slow,
  left-ward drift is not especially relevant in our moving frame
  (which focuses on the movement of $\Phi_t(\kappa(\gamma_0))$
  relative to that of $\Phi_t(\gamma_0)$). On successive laps
  around the torus, the orbit in question returns to the part of
  the torus shown in the figure, and the sequence of actions
  depicted in Fig.~\ref{quasiper snapshots} is repeated.  We
  regard this as geometric evidence of shear-induced chaos.

We have seen that in the phase-locked regime, the folding of the
limit cycle (when the time interval between kicks is sufficiently large) can
be deduced from the geometry of the strong stable foliation.  A
natural question is: in the quasi-periodic regime, are there
geometric clues in the unforced dynamics that will tell us
whether the system is predisposed to chaotic behavior when
forced?  Since folding occurs in finite time, we believe the
answer lies partially in what we call {\it finite-time stable
  manifolds}, a picture of which is shown in
Fig.~\ref{finite-time-wss}.  We first explain what these
manifolds are before discussing what they can --- and cannot ---
tell us.

Fix $t>0$.  At each $z \in {\mathbb T}^2$, let $V(z)$ be the
most contracted direction of the linear map $D\Phi_t(z)$ if it
is uniquely defined, {\em i.e.} if $v$ is a unit tangent vector
at $z$ in the direction $V(z)$, then $|D\Phi_t(z)v| \leq
|D\Phi_t(z)u|$ for all unit tangent vectors $u$ at $z$.  A
smooth curve is called a {\it time-$t$ stable manifold} if it is
tangent to $V$ at all points; these curves together form the
{\it time-$t$ stable foliation}.  In general, time-$t$ stable
manifolds are not necessarily defined everywhere; they vary with $t$,
and may
not stabilize as $t$ increases.  When ``real'' ({\it i.e.}
infinite-time) stable manifolds exist, time-$t$ stable
manifolds converge to them as $t \to \infty$.

The blue curve in Fig.~\ref{finite-time-wss} is an orbit segment
of $\Phi_t$.  The angles between this segment and the time-$5$
stable manifolds ({\em green curves}) reflect the presence of
shear. For example, if a kick sends points on the blue curve to
the right, then within $5$ units of time most points on the
kicked segment will lag behind their counterparts on the
original orbit segment --- except for the point with
$\theta_1=0$ at the time of the kick. Pinching certain points on
an orbit segment while having the rest slide back potentially
creates a scenario akin to that in Fig. 2; see
Sect.~\ref{sec1.2}.  One is also likely to find shear along the
black curve in Fig.~\ref{finite-time-wss}, a second orbit
segment of $\Phi_t$.  Whether or not the shear here is strong
enough to cause the formation of folds in 5 units of time cannot
be determined from the foliation alone; more detailed
information such as contraction rates are needed.  What
Fig.~\ref{finite-time-wss} tells us are the mechanism and the
shapes of the folds if they {\em do} form. Notice also that
shearing occurs in opposite directions along the blue and black
segments.  This brings us to a complication not present
previously: each orbit of $\Phi_t$ spends only a finite amount
of time near, say, the blue curve before switching to the region
near the black curve, and when it does so, it also switches the
direction of shear.  Finite-time stable foliations for system
(\ref{eq:2osc}) have also been computed for $t\in\{3,5\}$ and a
sample of $\afb \in (1.1, 1.6)$ (not shown).  They are
qualitatively similar to Fig.~\ref{finite-time-wss}, with most
of the leaves running in a northeasterly direction.

In summary, for $t$ not too large, time-$t$ stable foliations
generally do not change quickly with $t$ or with system
parameters.  They are good indicators of shear, but do not tell
us if there is {\it enough} shear for folds to form.  For the
system defined by (\ref{eq:2osc}), given that the finite-time
stable manifolds are nearly parallel to flowlines and the kick
map acts unevenly with respect to this foliation, we conclude
the presence of shear.  Fig.~\ref{quasiper snapshots} and
similar figures for other $\afb$ (not shown) confirm that
folding does indeed occur when the system is forced in the
near-periodic regime.

\heading{Computation of Lyapunov exponents}

To provide quantitative evidence of shear-induced chaos in the
situations discussed above, we compute $\Lambda_{\max}$.  Recall
that while periodic kicks followed by long relaxations provide a
simple setting to visualize folding, it is not expected to give
clean results for $\Lambda_{\max}$ because of the competition
between transient and sustained chaos (see Sect.~\ref{sec1.2}).
Continuous-time random forcing, on the other hand, produces
numerical results that are much easier to interpret.  

\paragraph{Study 4a: Stochastic Forcing.}
We consider system (\ref{eq:2osc}) with $\aff=1$ and $\afb \in
[1.1, 1.6]$.  The forcing is of the form $I(t)=a \cdot dB_t$
where $B_t$ is standard Brownian motion.

\paragraph{Study 4b: Periodic kicks.} 
The equation and parameters are as above, and the forcing is
given by $I(t) = A \cdot \sum_n \delta(t-nT)$.

\begin{figure}[h]
\begin{center}
\begin{tabular}{cc}
\includegraphics[bb=0in 0in 2.8in 2.5in]{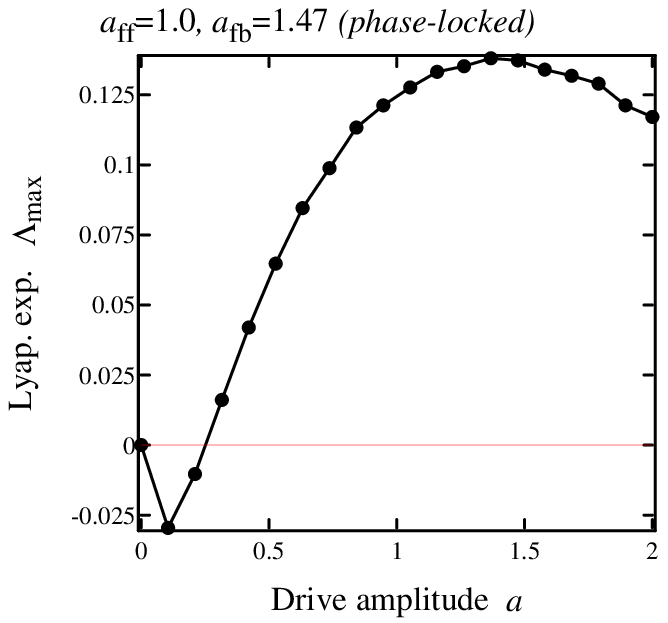}&
\includegraphics[bb=0in 0in 2.8in 2.5in]{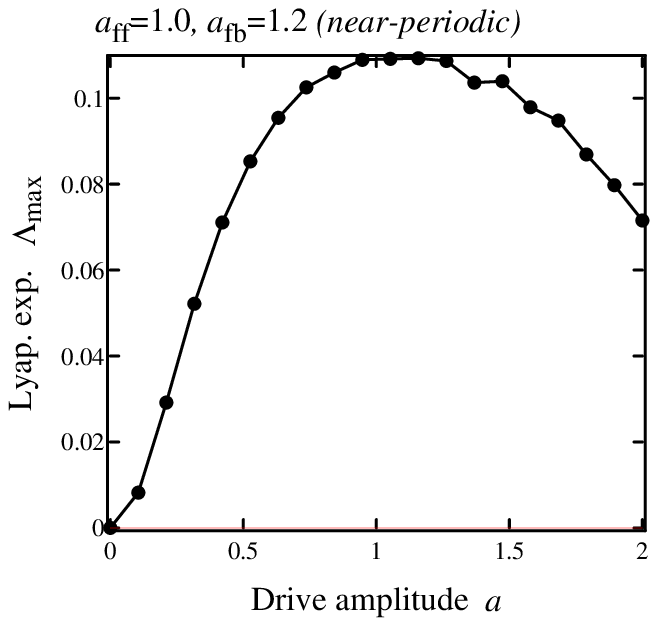}\\
(a)&(b)\\
\end{tabular}
\end{center}
\kaption{Lyapunov exponent of the system (\ref{eq:2osc})
  subjected to white noise forcing.  The parameters correspond
  to those in Figs.~\ref{cycle+wss} and \ref{quasiper
    snapshots}, respectively.}
\label{osc lyaps white noise}
\end{figure}

\begin{figure}[h!]
\begin{center}
\begin{tabular}{cc}
\includegraphics[bb=0in 0in 2.8in 2.5in]{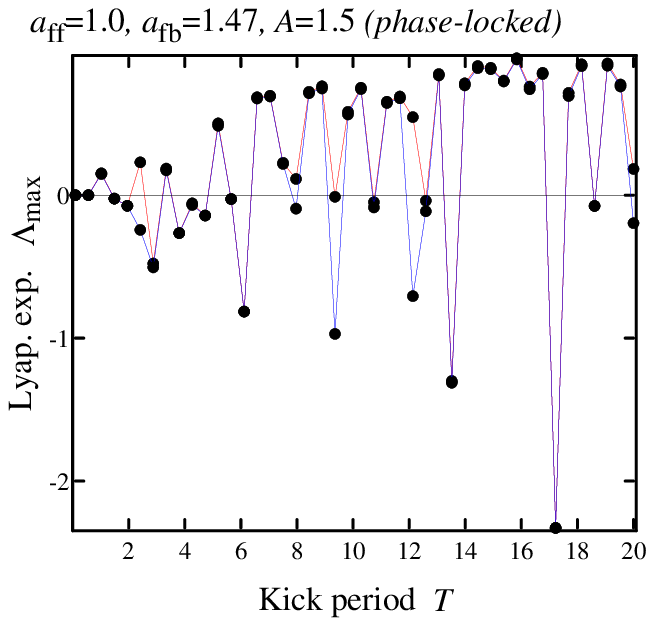}&
\includegraphics[bb=0in 0in 2.8in 2.5in]{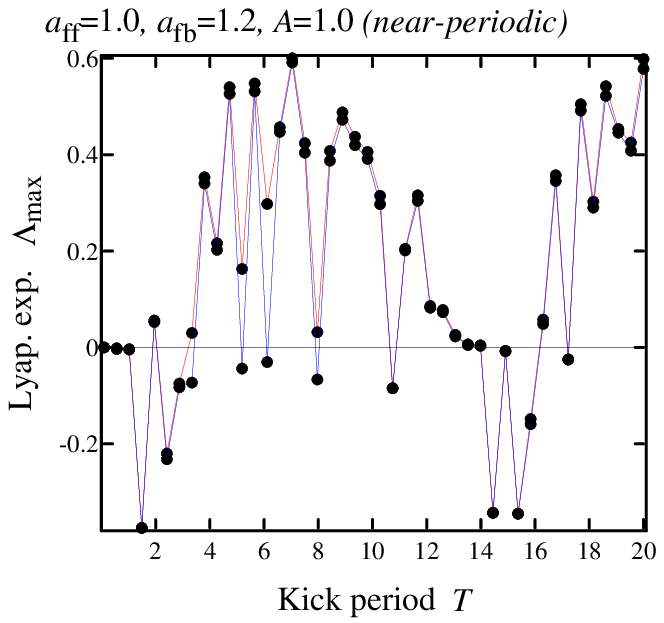}\\
(a)&(b)\\
\end{tabular}
\end{center}
\kaption{Lyapunov exponent of the system (\ref{eq:2osc})
  subjected to periodic kicks.  The parameters correspond to
  those in Figs.~\ref{cycle+wss} and \ref{quasiper snapshots},
  respectively.  As in Study 1, we show both upper and lower
  estimates of $\lmax$.}
\label{osc lyaps periodic kicks}
\end{figure}

\paragraph{Summary of Findings.}
{\it Positive $\Lambda_{\max}$ are found for stochastic forcing
  in the parameter interval studied, both before and after the
  onset of phase-locking at $\afb=\afb^*$. For periodic kicks with
  large enough $A$ and $T$, it appears that $\Lambda_{\max}$ is positive for a fraction
  of the forcing periods tested, but the results are hard to interpret
  due to the competition between
  transient and sustained chaos.}

\bigskip
\noindent
{\em Supporting Numerical Evidence.} Fig.~\ref{osc lyaps white
  noise} shows some results for stochastic forcing.  For
$\afb=1.47$, negative Lyapunov exponents are found for very
small amplitudes of forcing, while slightly stronger forcing
({\em e.g.} $a \approx 0.4$) is needed before $\Lambda_{\max}>0$
can be concluded with confidence.  In contrast, even fairly
small values of forcing seem to lead to $\lmax > 0$ when
$\afb=1.2$, {\em i.e.} in the near-periodic regime.  This may be
explained by the damping in the limit cycle case, especially for
larger $\afb$. Notice also that in this model large amplitudes
of forcing do not lead to larger $\Lambda_{\max}$.  This is due
to the fact that unlike the system in Studies 1--3, a very
strong forcing merely presses most of the phase space against
the circle $\theta_1=0$, which is not very productive from the
point of view of folding phase space.  Fig.~\ref{osc lyaps
  periodic kicks} shows plots of $\Lambda_{\max}$ for periodic
kicks.  Here, roughly 40\% of the kick periods $T$ for which
Lyapunov exponents were computed yield a positive exponent.
More generally, we find that $\lmax > 0$ for over 40\% of kick
intervals $T$ as $A$ varies over the range $[0.75,1.5]$. See
Simulation Details in Study 1.

%%%%%%%%%%%%%%%%%%%%%%%%%%%%
\section*{Conclusions}

Shear-induced chaos, by which we refer to the phenomenon of an
external force interacting with the shearing in a system to
produce stretches and folds, is found to occur for wide ranges
of parameters in forced oscillators and quasi-periodic systems.
Highlights of our results include:
\begin{itemize}

\item[(i)] For periodically kicked oscillators, positive
  Lyapunov exponents are observed under quite modest impositions
  on the unforced system and on the relaxation time between
  kicks (in contrast to existing rigorous results).  These
  regimes are, as expected, interspersed with those of transient
  chaos in parameter space.

\item[(ii)] Continuous-time stochastic forcing is shown to be
  equally effective in producing chaos.  The qualitative
  dependence on parameters is similar to that in deterministic
  forcing.  We find that suitably directed, degenerate white
  noise is considerably more effective than isotropic white
  noise (and additive noise will not work).  We have also found
  evidence for an approximate scaling law relating $\lmax$ to
  $\sigma$, $\lambda$, and $a$.  Other types of random forcing
  such as Poisson kicks are also studied and found to produce
  chaos.

\item[(iii)] The shear-induced stretching-and-folding mechanism
  can operate as well in quasi-periodic systems as it does in
  periodic systems, {\em i.e.} limit cycles are not a
  precondition for shear-induced chaos.  We demonstrate this
  through a pulse-coupled 2-oscillator system.  Chaos is induced
  under both periodic and white noise forcing, and a geometric
  explanation in terms of finite-time stable manifolds is
  proposed.

\end{itemize}
The conclusions in (i) and (ii) above are based on systematic
numerical studies of a linear shear flow model.  As this model
captures the essential features of typical oscillators, we
expect that our conclusions are valid for a wide range of other
models.  Our numerical results, particularly those on stochastic
forcing, point clearly to the possibility of a number of
(rigorous) theorems.


\begin{thebibliography}{42}

\bibitem{arnbook} L. Arnold, {\em Random Dynamical Systems},
  Springer-Verlag (1998)

\bibitem{arn} L. Arnold, N. Sri Namachchivaya,
  K. R. Schenk-Hopp\'e, ``Toward an understanding of stochastic
  Hopf bifurcation: a case study,'' {\em Int. J. Bifur. and
    Chaos} {\bf 6} (1996) {\em pp.} 1947--1975

\bibitem{aus} E. I. Auslender and G. N. Mil'shte{\u\i}n,
  ``Asymptotic expansions of the Liapunov index for linear
  stochastic systems with small noise,'' {\em
    J. Appl. Math. Mech.}  {\bf 46} (1982) {\em pp. 358--365}

\bibitem{bax94} P. H. Baxendale, ``A stochastic Hopf
  bifurcation,'' {\em Probab. Theory and Related Fields} (1994)
  {\em pp.} 581--616

\bibitem{bax03} P. H. Baxendale, ``Lyapunov exponents and
  stability for the stochastic Duffing-van der Pol oscillator,''
  {\em IUTAM Symposium on Nonlinear Stochastic Dynamics}, Kluwer
  (2003) {\em pp.} 125--135

\bibitem{bax04} P. H. Baxendale, ``Stochastic averaging and
  asymptotic behavior of the stochastic Duffing-van der Pol
  equation,'' {\em Stochastic Process. Appl.} {\bf 113}
  (2004) {\em pp.} 235--272

\bibitem{bax02} P. H. Baxendale, ``Lyapunov exponents and
  resonance for small periodic and random perturbations of a
  conservative linear system,'' {\em Stoch. Dyn.} {\bf 2} (2002)
  {\em pp.} 49--66

\bibitem{bax-gouk} P. H. Baxendale and L. Goukasian, ``Lyapunov
  exponents for small random perturbations of Hamiltonian
  systems,'' {\em Annals of Probability} {\bf 30} (2002) {\em
    pp.}  101--134

\bibitem{guck75} J. Guckenheimer, ``Isochrons and phaseless
  sets,'' {\em J. Theor. Biol.} {\bf 1} (1974) {\em pp.} 259--273

\bibitem{gh} J. Guckenheimer and P. Holmes, {\em Nonlinear
  Oscillations, Dynamical Systems, and Bifurcations of Vector
  Fields}, Springer-Verlag (1983)

\bibitem{gwy} J. Guckenheimer, M. Weschelberger, and
  L.-S. Young, ``Chaotic attractors of relaxation oscillators,''
  {\em Nonlinearity} {\bf 19} (2006) {\em pp.} 701--720

\bibitem{kifer} Yu. Kifer, {\em Ergodic Theory of Random
  Transformations}, Birkh\"{a}user (1986)

\bibitem{kunita} H. Kunita, {\em Stochastic Flows and Stochastic
  Differential Equations}, Cambridge University Press (1990)

\bibitem{prl} K. K. Lin, E. Shea-Brown, L.-S. Young, ``Reliable
  and unreliable dynamics in driven coupled oscillators,''
  preprint (2006); arXiv:nlin/0608021v1

\bibitem{reliability} K. K. Lin, E. Shea-Brown, L.-S. Young,
  ``Reliable and unreliable dynamics in coupled oscillator
  networks,'' in preparation

\bibitem{nsri} N. Sri Namachchivaya, ``The asymptotic stability
  of a weakly perturbed 2-dimensional non-Hamiltonian system",
  private communication

\bibitem{oksa+wang} A. Oksasoglu and Q. Wang, ``Strange
  attractors in periodically-kicked Chua's circuit,'' {\em
    Int. J. Bifur. Chaos} {\bf 16} (2005) {\em pp.} 83--98

\bibitem{pinsky} M. Pinsky and V. Wihstutz, ``Lyapunov exponents
  of nilpotent It\^o systems,'' {\em Stochastics} {\bf 25}
  (1998) {\em pp.} 43--57

\bibitem{sh} K. R. Schenk-Hopp\'e, ``Bifurcation scenarios of
  the noisy Duffing-van der Pol oscillator,'' {\em Nonlinear
    Dynamics} {\bf 11} (1996) {\em pp.} 255--274

\bibitem{taylor-holmes} D. Taylor and P. Holmes, ``Simple models
  for excitable and oscillatory neural networks,'' {\em
    J. Math. Biol.} {\bf 37} (1998) {\em pp.} 419--446

\bibitem{van der pol} B. van der Pol and J. van der Mark,
  ``Frequency demultiplication,''{\em Nature} {\bf 120} (1927)
  {\em pp.} 363--364

\bibitem{WY1} Q. Wang and L.-S. Young, ``Strange attractors with
  one direction of instability,'' {\em Comm. Math. Phys.}  {\bf
    218} (2001) {\em pp.} 1--97

\bibitem{WY2} Q. Wang and L.-S. Young, ``From invariant curves
  to strange attractors,'' {\em Comm. Math. Phys.}  {\bf 225}
  (2002) {\em pp.} 275--304

\bibitem{WY3} Q. Wang and L.-S. Young, ``Strange attractors
  in periodically-kicked limit cycles and Hopf bifurcations,''
  {\em Comm. Math. Phys.}  {\bf 240} (2003) {\em pp.} 509--529

\bibitem{WY4} Q. Wang and L.-S. Young, ``Toward a theory of rank
  one attractors,'' {\em Annals of Mathematics} (to appear)

\bibitem{winfree} A. Winfree, {\em The Geometry of Biological
  Time, Second Edition}, Springer-Verlag (2000)

\bibitem{zaslav} G. Zaslavsky, ``The simplest case of a strange
  attractor,'' {\em Physics Letters} {\bf 69A} (1978) {\em pp.}
  145--147

\end{thebibliography}
\end{document}